\documentclass[nonblindrev]{informs2017}
\OneAndAHalfSpacedXI % current default line spacing
\usepackage{endnotes}
\let\footnote=\endnote
\let\enotesize=\normalsize
\def\notesname{Endnotes}%
\def\makeenmark{\hbox to1.275em{\theenmark.\enskip\hss}}
\def\enoteformat{\rightskip0pt\leftskip0pt\parindent=1.275em
\leavevmode\llap{\makeenmark}}

\usepackage{multirow}
\usepackage{float}
\usepackage{epstopdf}
\usepackage{appendix}
\usepackage{bm}
\usepackage{dsfont}

\usepackage{graphicx}
\usepackage{subcaption}
\usepackage{caption}
\usepackage{enumitem}
\usepackage{array}

\usepackage{cleveref}

\newcommand{\bmt}[1]{\tilde{\bm{#1}}}
\newcommand{\bmh}[1]{\hat{\bm{#1}}}

\usepackage{tabularx}
\usepackage{pict2e}
\newcolumntype{L}[1]{>{\raggedright\arraybackslash}p{#1}}
\newcolumntype{C}[1]{>{\centering\arraybackslash}p{#1}}
\newcolumntype{R}[1]{>{\raggedleft\arraybackslash}p{#1}}
\usepackage{color}
\usepackage[usenames,dvipsnames]{xcolor}
\definecolor{strcolor}{rgb}{0.6, 0.2, 0.6}
\definecolor{commentcolor}{rgb}{0.3125, 0.5, 0.3125}
\definecolor{keycol}{rgb}{0, 0, 1}

\usepackage{listings}
\lstdefinelanguage{Julia}%
{morekeywords={abstract,break,case,catch,const,continue,do,else,elseif,%
		end,export,false,for,function,immutable,import,importall,if,in,%
		macro,module,otherwise,quote,return,switch,true,try,type,typealias,%
		using,while},%
sensitive=true,%
alsoother={\$},%
morecomment=[l]\#,%
morecomment=[n]{\#=}{=\#},%
morestring=[s]{"}{"},%
morestring=[m]{'}{'},%
}[keywords,comments,strings]%

\usepackage{natbib}
 \bibpunct[, ]{(}{)}{,}{a}{}{,}%
 \def\bibfont{\small}%
\TheoremsNumberedThrough     % Preferred (Theorem 1, Lemma 1, Theorem 2)
\ECRepeatTheorems

\EquationsNumberedThrough    % Default: (1), (2), ...
\begin{document}
%%%%%%%%%%%%%%%%

%\setcounter{page}{1}
%\DOI{10.1287/xxxx.0000.0000}%

\RUNAUTHOR{Chen, Kuhn, and Wiesemann}

\RUNTITLE{Approximations of Data-Driven Chance Constraints}

% \TITLE{Bundling Information Goods of Decreasing Value}
\TITLE{On Approximations of Data-Driven \\ Chance Constrained Programs over Wasserstein Balls}

\ARTICLEAUTHORS{%
\AUTHOR{Zhi Chen}
\AFF{Department of Management Sciences, College of Business, City University of Hong Kong, Kowloon Tong, Hong Kong, \\ zhi.chen@cityu.edu.hk}
\AUTHOR{Daniel Kuhn}
\AFF{Risk Analytics and Optimization Chair, \'{E}cole Polytechnique F\'{e}d\'{e}rale de Lausanne, Lausanne, Switzerland, \\ daniel.kuhn@epfl.ch}
\AUTHOR{Wolfram Wiesemann}
\AFF{Imperial College Business School, Imperial College London, London, United Kingdom, \\ ww@imperial.ac.uk}
}
\ABSTRACT{Distributionally robust chance constrained programs minimize a deterministic cost function subject to the satisfaction of one or more safety conditions with high probability, given that the probability distribution of the uncertain problem parameters affecting the safety condition(s) is only known to belong to some ambiguity set. We study three popular approximation schemes for distributionally robust chance constrained programs over Wasserstein balls, where the ambiguity set contains all probability distributions within a certain Wasserstein distance to a reference distribution. The first approximation replaces the chance constraint with a bound on the conditional value-at-risk, the second approximation decouples different safety conditions via Bonferroni's inequality, and the third approximation restricts the expected violation of the safety condition(s) so that the chance constraint is satisfied. We show that the conditional value-at-risk approximation can be characterized as a tight convex approximation, which complements earlier findings on classical (non-robust) chance constraints, and we offer a novel interpretation in terms of transportation savings. We also show that the three approximations can perform arbitrarily poorly in data-driven settings, and that they are generally incomparable with each other.
}%

%\KEYWORDS{deterministic inventory theory; infinite linear programming duality;
%  existence of optimal policies; semi-Markov decision process; cyclic schedule}


\HISTORY{\today}

\maketitle

%%%%%%%%%%%%%%%%%%%%%%%%%%%%%%%%%%%%%%%%%%%%%%%%%%%%%%%%%%%%%%%%%%%%%%
% Samples of sectioning (and labeling) in OPRE
% NOTE: (1) \section and \subsection do NOT end with a period
%       (2) \subsubsection and lower need end punctuation
%       (3) capitalization is as shown (title style).
%
%\section{Introduction.}\label{intro} %%1.
%\subsection{Duality and the Classical EOQ Problem.}\label{class-EOQ} %% 1.1.
%\subsection{Outline.}\label{outline1} %% 1.2.
%\subsubsection{Cyclic Schedules for the General Deterministic SMDP.}
%  \label{cyclic-schedules} %% 1.2.1
%\section{Problem Description.}\label{problemdescription} %% 2.

% Text of your paper here

\section{Introduction}\label{sec:introduction}

In this paper we study data-driven distributionally robust chance constrained programs of the form
\begin{equation}\label{prob:cc general}
\begin{array}{cll}
\displaystyle \min_{\bm{x} \in \mathcal{X}} &~\bm{c}^\top\bm{x} \\
{\rm s.t.} &~\displaystyle \mathbb{P}[\bmt{\xi} \in \mathcal{S}(\bm{x})] \geq 1-\varepsilon &~\forall \mathbb{P} \in \mathcal{F}(\theta).
\end{array}
\end{equation}
The goal is to find a decision $\bm{x}$ from within a compact polyhedron $\mathcal{X} \subseteq \mathbb{R}^L$ that minimizes a linear cost function $\bm{c}^\top\bm{x}$ and ensures that the exogenous random vector $\bmt{\xi}$ falls within a decision-dependent safety set $\mathcal{S}(\bm{x}) \subseteq \mathbb{R}^K$ with high probability $1-\varepsilon$ under every distribution $\mathbb{P}$ that resides in the Wasserstein ball $\mathcal{F}(\theta)$ of radius $\theta\ge 0$:
\begin{equation*}\label{set:Wasserstein}
\mathcal{F}(\theta) = \{\mathbb{P} \in \mathcal{P}(\mathbb{R}^K) \mid
d_{\rm W}(\mathbb{P}, \hat{\mathbb{P}}) \leq \theta\}.
\end{equation*}
Here, $\hat{\mathbb{P}} = \frac{1}{N}\sum_{i = 1}^N \delta_{\bmh{\xi}_i}$ is the empirical distribution over $N$ historical samples $\{\bmh{\xi}_i\}_{i \in [N]}$ of $\bmt{\xi}$, and the (type-1) Wasserstein distance $d_{\rm W}(\mathbb{P}_1,\mathbb{P}_2)$ between two distributions $\bmt{\xi}_1 \sim \mathbb{P}_1$ and $\bmt{\xi}_2 \sim \mathbb{P}_2$ on $\mathbb{R}^K$, equipped with a general norm $\|\cdot\|$, is defined as
$$
\begin{array}{rcl}
d_{\rm W}(\mathbb{P}_1,\mathbb{P}_2) \; = \; &\displaystyle \inf_{\mathbb{P} \in \mathcal{P}(\mathbb{P}_1, \mathbb{P}_2)}
& \mathbb{E}_{\mathbb{P}}[\|\bmt{\xi}_1 - \bmt{\xi}_2\|],
\end{array}
$$
where $\mathcal{P}(\mathbb{P}_1, \mathbb{P}_2)$ is the set of all joint distributions on $\mathbb{R}^K\times \mathbb{R}^K$ with marginals $\mathbb{P}_1$ and $\mathbb{P}_2$. Problem~\eqref{prob:cc general} generalizes both individual chance constrained programs, where $\mathcal{S}(\bm{x}) = \{ \bm{\xi} \in \mathbb{R}^K \mid (\bm{A} \bm{\xi} + \bm{a})^\top \bm{x} < \bm{b}^\top \bm{\xi} + b_0 \}$ for $\bm{A} \in \mathbb{R}^{L \times K}$, $\bm{a} \in \mathbb{R}^L$, $\bm{b} \in \mathbb{R}^K$ and $b_0 \in \mathbb{R}$, and joint chance constrained programs with right-hand side uncertainty, where $\mathcal{S}(\bm{x}) = \{\bm{\xi} \in \mathbb{R}^K \mid \bm{a}_m^\top \bm{x} < \bm{b}_m^\top \bm{\xi} + b_{m0} ~\forall m \in [M] \}$ for $\bm{a}_m \in \mathbb{R}^L$, $\bm{b}_m \in \mathbb{R}^K$ and $b_{m0} \in \mathbb{R}$, $m \in [M]$.

It has been shown that a fixed decision $\bm{x}$ satisfies the ambiguous chance constraint in~\eqref{prob:cc general} if and only if the partial sum of the $\varepsilon N$ smallest transportation distances to the unsafe set $\bar{\mathcal{S}}(\bm{x}) = \mathbb{R}^K \setminus \mathcal{S}(\bm{x})$, multiplied by the mass $1/N$ of a training sample, exceeds~$\theta$.

\begin{theorem}[\textnormal{\citet{Chen_Kuhn_Wiesemann_2022}}]\label{thm:cc equivalent}
For any fixed decision $\bm{x} \in \mathcal{X}$, the ambiguous chance constraint in~\eqref{prob:cc general} is satisfied if and only if
\begin{equation*}\label{equivalence:theta positive}
\dfrac{1}{N}\sum_{i = 1}^{\varepsilon N} \mathbf{dist}(\bmh{\xi}_{\pi_i(\bm{x})}, \bar{\mathcal{S}}(\bm{x})) \ge \theta.
\end{equation*}
Here, $\bm{\pi}(\bm{x}) : [N] \rightarrow [N]$ is a decision-dependent permutation that orders the training samples $\{\bmh{\xi}_i\}_{i \in [N]}$ in order of non-decreasing distance to the unsafe set $\bar{\mathcal{S}}(\bm{x})$, and the distance with respect to a norm $\|\cdot\|$ is defined as $\mathbf{dist}(\bmh{\xi}_i, \bar{\mathcal{S}}(\bm{x})) = \min\{\|\bm{\xi}-\bmh{\xi}_i\| \mid \bm{\xi} \in \bar{\mathcal{S}}(\bm{x})\}$.
\end{theorem}

(The sum in Theorem~\ref{thm:cc equivalent} is defined even if $\varepsilon N \notin \mathbb{N}$; please refer to the notation at the end of this section.) Theorem~\ref{thm:cc equivalent} allows us to reformulate individual and joint chance constrained programs as deterministic mixed-integer conic programs \citep{Chen_Kuhn_Wiesemann_2022, xie2021distributionally}.

\begin{proposition}[\textnormal{\citet{xie2021distributionally, Chen_Kuhn_Wiesemann_2022}}]\label{prop:individual cc}
For the safety set $\mathcal{S}(\bm{x}) = \{\bm{\xi} \in \mathbb{R}^K \mid (\bm{A}\bm{\xi} + \bm{a})^\top \bm{x} < \bm{b}^\top\bm{\xi} + b_0\}$, where $\bm{A}^\top\bm{x} \ne \bm{b}$ for all $\bm{x} \in \mathcal{X}$, problem~\eqref{prob:cc general} is equivalent to the mixed-integer conic program
\begin{equation*}\label{prob:individual cc reformulation linearization}
\begin{array}{rcll}
Z^\star_{\rm ICC} =& \displaystyle \min_{\bm{q}, \bm{s}, t, \bm{x}} & \bm{c}^\top\bm{x} \\
&{\rm s.t.} & \varepsilon N t - \mathbf{e}^\top\bm{s} \geq \theta N \|\bm{b} - \bm{A}^\top\bm{x}\|_* \\
&& (\bm{b} - \bm{A}^\top\bm{x})^\top\bmh{\xi}_i + b_0 - \bm{a}^\top\bm{x} + {\rm M} q_i \geq t - s_i &~\forall i \in [N] \\
&& {\rm M} (1 - q_i) \geq t - s_i &~\forall i \in [N] \\
&& \bm{q} \in \{0,1\}^N, ~\bm{s} \geq \bm{0}, ~\bm{x} \in \mathcal{X},
\end{array}
\end{equation*}
where ${\rm M}$ is a suitably large (but finite) positive constant.
\end{proposition}

The condition that $\bm{A}^\top\bm{x} \ne \bm{b}$ for all $\bm{x} \in \mathcal{X}$ in Proposition~\ref{prop:individual cc} is non-restrictive. In fact, the weaker condition that $\bm{A}^\top\bm{x}^\star \ne \bm{b}$ for any \emph{optimal} solution $\bm{x}^\star \in \mathcal{X}$ is sufficient, and if an optimal solution $\bm{x}^\star$ satisfies $\bm{A}^\top\bm{x} = \bm{b}$, then an alternative optimal solution $\bm{x}'$ satisfying $\bm{A}^\top\bm{x}' \ne \bm{b}$ can be identified from the solution of auxiliary optimization problems.
We refer to \citet{Chen_Kuhn_Wiesemann_2022} for the details.

\begin{proposition}[\textnormal{\citet{xie2021distributionally, Chen_Kuhn_Wiesemann_2022}}]\label{prop:joint cc}
For the safety set $\mathcal{S}(\bm{x}) = \{\bm{\xi} \in \mathbb{R}^K \mid \bm{a}^\top_m \bm{x} < \bm{b}^\top_m\bm{\xi} + b_{m0} ~\forall m \in [M]\}$, where $\bm{b}_m \ne \bm{0}$ for all $m \in [M]$, problem~\eqref{prob:cc general} is equivalent to the mixed-integer conic program
\begin{equation*}
\label{prob:joint cc reformulation M linearization}
\begin{array}{rcll}
Z^\star_{\rm JCC} =& \displaystyle \min_{\bm{q}, \bm{s}, t, \bm{x}} & \bm{c}^\top\bm{x} \\
&{\rm s.t.} & \varepsilon N t - \mathbf{e}^\top\bm{s} \geq \theta N \\
&& \dfrac{\bm{b}^\top_m\bmh{\xi}_i + b_{m0} - \bm{a}^\top_m\bm{x}}{\|\bm{b}_m\|_*} + {\rm M} q_i \geq t - s_i &~\forall m \in [M], ~i \in [N] \\
&& {\rm M} (1 - q_i) \geq t - s_i &~\forall i \in [N] \\
&& \bm{q} \in \{0,1\}^N,~ \bm{s} \geq \bm{0}, ~\bm{x} \in \mathcal{X},
\end{array}
\end{equation*}
where ${\rm M}$ is a suitably large (but finite) positive constant.
\end{proposition}

The assumption that $\bm{b}_m \ne \bm{0}$ for all $m \in [M]$ in Proposition~\ref{prop:joint cc} is non-restrictive since any safety condition with $\bm{b}_m = \bm{0}$ is deterministic and can thus be absorbed in the definition of the set $\mathcal{X}$.

In recent years, the mixed-integer conic programming reformulations developed in Propositions~\ref{prop:individual cc} and~\ref{prop:joint cc} have been strengthened so as to scale better to larger problem sizes \citep{ho2020strong, ho2021distributionally, ji2021data, shen2022chance}. Nevertheless, similar to classical chance constraints as well as distributionally robust chance constraints over moment ambiguity sets, exact reformulations of the distributionally robust chance constrained program~\eqref{prob:cc general} quickly become computationally prohibitive for large problems. As a result, there has been significant interest in safe (\emph{i.e.}, conservative) tractable approximations to problem~\eqref{prob:cc general} that scale gracefully with problem size.

Distributionally robust chance constrained programs are most commonly approximated by the Bonferroni approximation or the worst-case conditional value-at-risk (CVaR) approximation. The quality of the Bonferroni approximation crucially depends on the choice of the associated Bonferroni weights. While \cite{xie2019optimized} show that these Bonferroni weights can be optimized efficiently under specific conditions, \cite{Chen_Sim_Sun_Teo_2010} show that the quality of the Bonferroni approximation can be poor even if the Bonferroni weights are chosen optimally. \cite{Chen_Sim_Sun_Teo_2010} also show that the worst-case CVaR approximation can outperform the Bonferroni approximation with optimally chosen Bonferroni weights for Chebyshev (\emph{i.e.}, second-order moment) ambiguity sets, provided that certain scaling factors in the worst-case CVaR approximation are selected judiciously. \cite{zymler2013distributionally} show that the worst-case CVaR approximation is indeed exact for distributionally robust chance constrained programs over Chebyshev ambiguity sets if the scaling factors are selected optimally. This result has been extended to non-linear safety conditions by \cite{yang2016distributionally}. Selecting the scaling factors optimally, however, amounts to solving a non-convex optimization problem. More recently, \cite{ahmed2017nonanticipative} and \cite{jiang2022also} have proposed the ALSO-X approximation, which traces the Pareto efficient solutions in terms of the objective value of problem~\eqref{prob:cc general} and the expected violation of the safety condition(s) to determine a feasible (but typically suboptimal) solution to~\eqref{prob:cc general}. It has been shown that the ALSO-X approximation can outperform the worst-case CVaR approximation in chance constrained programs with known distributions as well as type-$\infty$ Wasserstein ambiguity sets. For further information, we refer the reader to the surveys by \cite{Ben-tal_Nemirovski_book}, \cite{nemirovski2012safe} and \cite{Hanasusanto_Roitch_Kuhn_Wiesemann_2015}. 

This paper complements the literature with the following three contributions.
\begin{enumerate}[leftmargin=*, labelindent=16pt]
\item We show that the CVaR offers a tight convex approximation to certain disjunctive constraints underlying the reformulations of Propositions~\ref{prop:individual cc} and~\ref{prop:joint cc}. This provides a theoretical justification for the popularity of this approximation scheme. 
\item We show that the CVaR approximation admits an intuitive interpretation: it accounts for transportation \emph{savings} alongside the transportation costs in the computation of the Wasserstein distance. This complements existing interpretations of chance constraints over Wasserstein balls based on the CVaR (see, \emph{e.g.}, \citealt{xie2021distributionally, ho2022adversarial}).
\item We show that the CVaR, Bonferroni and ALSO-X approximations may result in solutions that are severely suboptimal, and that these three approximation schemes are generally incomparable with each other. This contrasts earlier results for moment ambiguity sets, where the CVaR approximation is known to outperform the Bonferroni approximation under optimally chosen scaling factors, and for type-$\infty$ Wasserstein balls, where the ALSO-X approximation has been shown to be tighter than the CVaR approximation.
\end{enumerate}

\vspace{5mm}
\noindent \textbf{Notation.}
Boldface uppercase (resp., lowercase) letters denote matrices (resp., vectors). Special vectors of an appropriate dimension include $\bm{0}$ and $\mathbf{e}$, which represent the zero vector and the vector of all ones, respectively. We let $ [N] = \left\{1,2,\ldots,N\right\} $ denote the set of positive integers up to $ N $, and $\|\cdot\|_*$ represents the dual norm of a general norm $\|\cdot\|$. Given a possibly fractional number $\ell\in [0,N]$, the partial sum of the first $\ell$ values in $\{k_i\}_{i \in [N]}$ is defined as $\sum_{i = 1}^{\ell} k_i = \sum_{i = 1}^{\lfloor \ell \rfloor} k_i + (\ell - \lfloor \ell \rfloor) k_{\lfloor \ell \rfloor + 1}$. For $x \in \mathbb{R}$, we define $(x)^+ = \max \{ x, \, 0 \}$. Finally, we denote random vectors by tilde signs ({\em e.g.}, $\bmt{\xi}$) and their realizations by the same symbols without tildes ({\em e.g.}, $\bm{\xi}$).

\section{CVaR Approximation}

We propose a systematic approach to constructing safe convex approximations for problem~\eqref{prob:cc general} under individual and joint chance constraints in Sections~\ref{sec:cvar_individual_cc} and~\ref{sec:cvar_joint_cc}, respectively. We show that the CVaR approximation to the respective chance constraint represents the best among a low-parametric class of approximations, and we elucidate how the CVaR approximation can be interpreted as assigning transportation \emph{savings} to certain atoms of the empirical distribution.

\subsection{Individual Chance Constraints}\label{sec:cvar_individual_cc}

Consider an instance of problem~\eqref{prob:cc general} with an individual chance constraint corresponding to the safety set $\mathcal{S}(\bm{x}) = \{\bm{\xi} \in \mathbb{R}^K \mid (\bm{A}\bm{\xi} + \bm{a})^\top \bm{x} < \bm{b}^\top\bm{\xi} + b_0 \}$. As discussed in Section~\ref{sec:introduction}, we can assume that $\bm{A}^\top\bm{x} \ne \bm{b}$ for all $\bm{x} \in \mathcal{X}$. By \Cref{prop:individual cc}, the distributionally robust chance constrained program~\eqref{prob:cc general} is thus equivalent to the deterministic optimization problem 
\begin{equation}
\label{prob:individual cc reformulation abstract}
Z^\star_{\rm ICC}= \min_{(\bm{x},\bm{s}, t)\in \mathcal{C}_{\rm ICC}} \bm c^\top \bm x,
\end{equation}
whose feasible region is given by
\begin{equation*}
\mathcal{C}_{\rm ICC} = \left\{(\bm{x},\bm{s}, t) \in \mathcal{X} \times \mathbb{R}^N_+ \times \mathbb{R} ~\left|~
\begin{array}{l}
\varepsilon N t - \mathbf{e}^\top\bm{s} \geq \theta N \|\bm{b} - \bm{A}^\top\bm{x}\|_* \\
((\bm{b} - \bm{A}^\top\bm{x})^\top\bmh{\xi}_i + b_0 - \bm{a}^\top\bm{x})^+ \geq t - s_i \quad \forall i \in [N] 
\end{array}
\right.\right\}.
\end{equation*}	
Since $\mathcal{C}_{\rm ICC}$ is non-convex, it is natural to replace it with tractable conservative (inner) approximations. We next show that any convex inner approximation of  $\mathcal{C}_{\rm ICC}$ is dominated, in the sense of set inclusion, by a convex set of the form
\begin{equation*}
\mathcal{C}_{\rm ICC}(\bm{\kappa}) = \left\{(\bm{x},\bm{s}, t) \in \mathcal{X} \times \mathbb{R}^N_+ \times \mathbb{R} ~\left|~
\begin{array}{l}
\varepsilon N t - \mathbf{e}^\top\bm{s} \geq \theta N \|\bm{b} - \bm{A}^\top\bm{x}\|_* \\
\kappa_i ((\bm{b} - \bm{A}^\top\bm{x})^\top\bmh{\xi}_i + b_0 - \bm{a}^\top\bm{x}) \geq t - s_i \quad \forall i \in [N] 
\end{array}
\right.\right\}
\end{equation*}	
parameterized by a vector of slope parameters $\bm \kappa \in [0,1]^N$.

\begin{proposition}\label{prop:dominant convex approximation_1}
For any convex set $\mathcal{W} \subseteq \mathcal{C}_{\rm ICC}$, there exists $\bm{\kappa} \in[0,1]^N$ with $\mathcal{W} \subseteq \mathcal{C}_{\rm ICC} (\bm{\kappa}) \subseteq \mathcal{C}_{\rm ICC}$.
\end{proposition}

\begin{figure}[tb]
\begin{subfigure}{.5\textwidth}
\begin{center}
\includegraphics[width=0.75\linewidth]{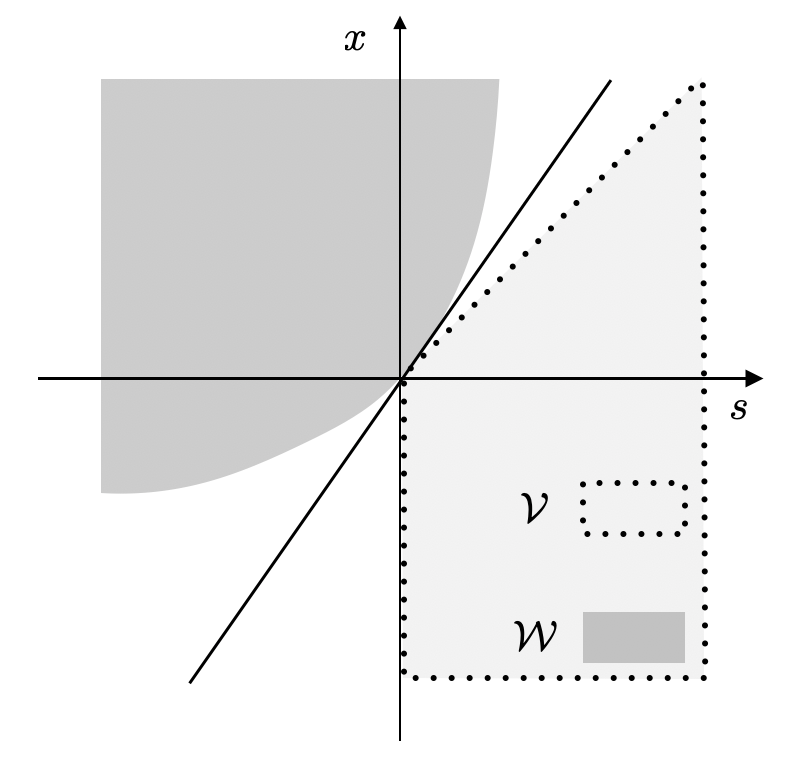}
\end{center}
\end{subfigure}%
\begin{subfigure}{.5\textwidth}
\begin{center}
\includegraphics[width=0.75\linewidth]{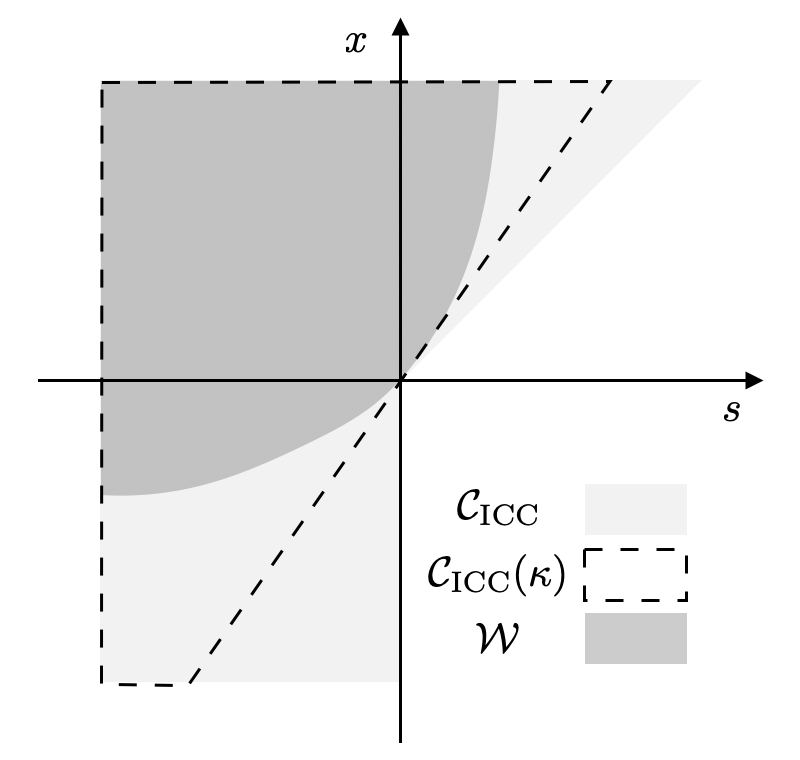}
\end{center}
\end{subfigure}
\vspace{0.2cm}
	
\caption{{\textnormal{The complement of $\mathcal{C}_{\rm ICC} = \{(x,s) \in \mathbb{R}^2 \mid x^+ \geq s\}$, denoted by $\mathcal{V} = \mathbb{R}^2 \setminus \mathcal{C}_{\rm ICC} = \{(x,s) \in \mathbb{R}^2 \mid x^+ < s\}$, is convex and intersection-free with any convex inner approximation $\mathcal{W}$ of $\mathcal{C}_{\rm ICC}$. Thus, $\mathcal{V}$ and $\mathcal{W}$ admit a separating hyperplane. The separating hyperplane can always be chosen such that it passes through the origin, that is, it satisfies $\kappa x = s$ for some $\kappa \in [0,1]$. The separating hyperplane then yields a convex set $\mathcal{C}_{\rm ICC}(\kappa) = \{(x,s) \in \mathbb{R}^2 \mid \kappa x \geq s\}$ that weakly dominates $\mathcal{W}$.}}} \label{fig:inner approximation}
\end{figure}

\Cref{prop:dominant convex approximation_1} implies that amongst all convex conservative approximations to problem~\eqref{prob:individual cc reformulation abstract} it is sufficient to focus on those that are induced by a feasible set of the form $\mathcal{C}_{\rm ICC}(\bm \kappa)$ for some $\bm{\kappa} \in[0,1]^N$. The intuition behind \Cref{prop:dominant convex approximation_1} is illustrated in Figure~\ref{fig:inner approximation}. Thus, it is sufficient to focus on the family of approximate problems of the form
\begin{equation}
\label{prob:linear approximation individual cc_1}
Z^\star_{\rm ICC}(\bm{\kappa}) =\min_{(\bm{x},\bm{s}, t)\in \mathcal{C}_{\rm ICC}(\bm \kappa)} \bm c^\top \bm x
\end{equation}
parameterized by $\bm{\kappa} \in[0,1]^N$. The following proposition asserts that the best approximation within this family is exact.

\begin{proposition}\label{prop:exact dominant convex approximation_1}
We have $Z^\star_{\rm ICC} = \min_{\bm\kappa \in [0,1]^N} Z^\star_{\rm ICC}(\bm{\kappa})$.
\end{proposition}

\Cref{prop:exact dominant convex approximation_1} implies that the family~\eqref{prob:linear approximation individual cc_1} of tractable upper bounding problems contains an instance $\bm \kappa^\star\in \arg\min_{\bm \kappa\in [0,1]^N} Z^\star_{\rm ICC}(\bm{\kappa})$ that recovers an optimal solution to the ambiguous chance constrained program~\eqref{prob:individual cc reformulation abstract}, which is known to be NP-hard \citep[Theorem~12]{xie2020bicriteria}. We may thus conclude that computing $\bm\kappa^\star$ is also NP-hard. The complexity of computing the best upper bound within the family~\eqref{prob:linear approximation individual cc_1} can be reduced by restricting attention to uniform slope parameters of the form $\bm{\kappa} = \kappa\mathbf{e}$ for some $\kappa \in [0,1]$. Within this subset, the choice~$\bm \kappa=\mathbf{e}$ turns out to be optimal.

\begin{proposition}\label{prop:cvar best_1}
We have $\min_{\kappa \in [0,1]} Z^\star_{\rm ICC}(\kappa\mathbf{e})= Z^\star_{\rm ICC}(\mathbf{e})$.
\end{proposition}

Next, we demonstrate that the approximate problem~\eqref{prob:linear approximation individual cc_1} corresponding to $\bm \kappa=\mathbf{e}$ can also be obtained by approximating the worst-case chance constraint in~\eqref{prob:cc general} with a worst-case CVaR constraint. To see this, note first that 
\begin{eqnarray*}
\mathbb{P}[\bmt{\xi} \in \mathcal{S}(\bm{x})] \ge 1-\varepsilon ~& \iff &~ \mathbb{P}[(\bm{A}\bmt{\xi} + \bm{a})^\top \bm{x} \ge \bm{b}^\top\bmt{\xi} + b_0]\le\varepsilon \\
~&\iff &~ \mathbb{P}\text{-VaR}_{\varepsilon}(\bm{a}^\top\bm{x} - b_0 + (\bm{A}^\top\bm{x}-\bm{b})^\top\bmt{\xi}) \leq 0 \\
~&\Longleftarrow &~ \mathbb{P}\text{-CVaR}_{\varepsilon}(\bm{a}^\top\bm{x} - b_0 + (\bm{A}^\top\bm{x}-\bm{b})^\top\bmt{\xi}) \leq 0
\end{eqnarray*}
for any $\mathbb P\in  \mathcal{F}(\theta)$, where the $\varepsilon$-value-at-risk (VaR) and the $\varepsilon$-CVaR of a measurable loss function $\ell (\bm{\xi})$ are defined as $\mathbb{P}\text{-VaR}_{\varepsilon} (\ell (\bm{\xi})) = \inf\{\gamma \in \mathbb{R} \mid \mathbb{P}[\gamma \leq \ell(\bmt{\xi})] \leq \varepsilon\}$ and $\mathbb{P}\text{-CVaR}_{\varepsilon}(\ell(\bmt{\xi})) = \inf\{\tau + \mathbb{E}_{\mathbb{P}}[(\ell(\bmt{\xi}) - \tau)^+]/\varepsilon \mid \tau \in \mathbb{R}\}$, respectively. The first equivalence above follows from the definition of the safety set, the second equivalence holds due to the definition of the VaR, and the last implication exploits the fact that the CVaR upper bounds the VaR. Thus, the worst-case CVaR constrained program
\begin{equation}
\label{eq:wc-cvar-program}
Z^\star_{\rm CVaR}=\left\{ \begin{array}{cl}
\displaystyle \min_{\bm{x} \in \mathcal{X}} & \bm{c}^\top\bm{x} \\
{\rm s.t.} & \displaystyle \sup_{\mathbb{P} \in \mathcal{F}(\theta)} \mathbb{P}\text{-CVaR}_{\varepsilon}(\bm{a}^\top\bm{x} - b_0 + (\bm{A}^\top\bm{x}-\bm{b})^\top\bmt{\xi}) \leq 0
\end{array}\right.
\end{equation}
constitutes a conservative approximation for the worst-case chance constrained program~\eqref{prob:cc general}, that is, $Z^\star_{\rm ICC} \le Z^\star_{\rm CVaR}$. We are now ready to prove that $Z^\star_{\rm CVaR}= Z^\star_{\rm ICC}(\mathbf{e})$.

\begin{proposition}\label{prop:individual cc worst-case CVaR approximation_1}
We have $Z^\star_{\rm CVaR}= Z^\star_{\rm ICC}(\mathbf{e})$.
\end{proposition}

\begin{remark}
Using similar arguments as in \Cref{prop:individual cc worst-case CVaR approximation_1}, one can show that problem~\eqref{prob:linear approximation individual cc_1} with $\bm \kappa=\kappa\mathbf{e}$ for any $\kappa\in (0,1]$ is equivalent to a worst-case CVaR constrained program of the form~\eqref{eq:wc-cvar-program}, where the Wasserstein radius $\theta$ is inflated to $\theta/\kappa$. This observation reconfirms that $\kappa=1$ is the least conservative choice amongst all uniform slope parameters in~\eqref{prob:linear approximation individual cc_1}; see \Cref{prop:cvar best_1}.
\end{remark}

The intimate links between the worst-case CVaR approximation~\eqref{eq:wc-cvar-program} and the worst-case chance constrained program~\eqref{prob:cc general} can also be studied through the lens of \Cref{thm:cc equivalent}. To this end, recall that the ambiguous chance constraint in problem~\eqref{prob:cc general} is equivalent to the deterministic constraint
\begin{equation*}
\dfrac{1}{N}\sum_{i = 1}^{\varepsilon N} \mathbf{dist}(\bmh{\xi}_{\pi_i(\bm{x})}, \bar{\mathcal{S}}(\bm{x})) \ge \theta.
\end{equation*}
We define the \emph{signed distance} between a point $\bm{\xi}$ and a closed set $\mathcal{C}$ as $\mathbf{sgn\textbf{-}dist} (\bm{\xi}, \mathcal{C}) = \mathbf{dist} (\bm{\xi}, \mathcal{C})$ if $\bm{\xi} \notin \mathcal{C}$ and $\mathbf{sgn\textbf{-}dist} (\bm{\xi}, \mathcal{C}) = -\mathbf{dist} (\bm{\xi}, \text{bd}(\mathcal{C}))$ otherwise. Here, $\text{bd}(\mathcal{C})$ denotes the boundary of the set $\mathcal{C}$. We then obtain the following result.

\begin{proposition}\label{prop:CVaR equivalence_1}
For any fixed decision $\bm{x} \in \mathcal{X}$, we have
$$
\sup_{\mathbb{P} \in \mathcal{F}(\theta)} \mathbb{P}\text{\emph{-CVaR}}_{\varepsilon}(\bm{a}^\top\bm{x} - b_0 + (\bm{A}^\top\bm{x}-\bm{b})^\top\bmt{\xi}) \leq 0
~\Longleftrightarrow~ 
\dfrac{1}{N}\sum_{i = 1}^{\varepsilon N} \mathbf{sgn\textbf{-}dist} (\bmh{\xi}_{\pi_i(\bm{x})}, \bar{\mathcal{S}}(\bm{x})) \ge \theta,
$$
where $\bm{\pi}(\bm{x})$ permutes the data points $\bmh{\xi}_i$ into ascending order of their signed distances to $\bar{\mathcal{S}}(\bm{x})$.
\end{proposition}

\Cref{thm:cc equivalent} and \Cref{prop:CVaR equivalence_1} show that both the ambiguous chance constraint in problem~\eqref{prob:cc general} and its worst-case CVaR approximation~\eqref{eq:wc-cvar-program} impose lower bounds on the costs of moving a fraction $\varepsilon$ of the training samples to the unsafe set. Moreover, since $\mathbf{sgn\textbf{-}dist} (\bmh{\xi}_i, \bar{\mathcal{S}} (\bm{x})) \leq \mathbf{dist} (\bmh{\xi}_i, \bar{\mathcal{S}} (\bm{x}))$ by construction, the worst-case CVaR constraint conservatively approximates the ambiguous chance constraint. In fact, we have $\mathbf{sgn\textbf{-}dist} (\bmh{\xi}_i, \bar{\mathcal{S}} (\bm{x})) = \mathbf{dist} (\bmh{\xi}_i, \bar{\mathcal{S}} (\bm{x}))$ for safe scenarios $\bmh{\xi}_i \in \mathcal{S} (\bm{x})$, whereas $\mathbf{sgn\textbf{-}dist} (\bmh{\xi}_j, \bar{\mathcal{S}} (\bm{x})) < 0$ even though $\mathbf{dist} (\bmh{\xi}_j, \bar{\mathcal{S}} (\bm{x})) = 0$ for (strictly) unsafe scenarios $\bmh{\xi}_j \in \text{int}(\bar{\mathcal{S}} (\bm{x}))$. In other words, the worst-case CVaR approximation~\eqref{eq:wc-cvar-program} assigns \emph{fictitious transportation savings} to training samples that are contained in the unsafe set. This leads to the following insight.

\begin{corollary}\label{prop:individual cc CVaR is exact_1}
The worst-case CVaR approximation is exact, that is, $Z^\star_{\rm CVaR}= Z^\star_{\rm ICC}$, under either of the following conditions.
\begin{itemize}
\item[(i)] We have $\bmh{\xi}_i \in \mathcal{S}(\bm{x}^\star)$ for all $i \in [N]$, where $\bm{x}^\star$ is optimal in~\eqref{prob:cc general}.
\item[(ii)] We have $\varepsilon \leq 1/N$. 
\end{itemize}
\end{corollary}

\begin{remark}
Propositions~\ref{prop:exact dominant convex approximation_1},~\ref{prop:cvar best_1} and~\ref{prop:individual cc worst-case CVaR approximation_1} show that the worst-case CVaR approximation, despite being a best-in-class approximation among the inner linearizations $\mathcal{C}_{\textrm{ICC}}$ with uniform slopes, is in general not exact, that is,
\[
Z^\star_{\rm ICC}
\;\; = \;\;
\min_{\bm{\kappa} \in [0,1]^N} Z^\star_{\rm ICC}(\bm{\kappa})
\;\; \leq \;\;
\min_{\kappa \in [0,1]}Z^\star_{\rm ICC}(\kappa\mathbf{e})
\;\; = \;\;
Z^\star_{\rm ICC}(\mathbf{e})
\;\; = \;\;
Z^\star_{\rm CVaR}.
\]
\Cref{prop:individual cc CVaR is exact_1} identifies sufficient conditions for the worst-case CVaR approximation to be exact.
\end{remark}

\subsection{Joint Chance Constraints with Right-Hand Side Uncertainty}\label{sec:cvar_joint_cc}

Consider now an instance of problem~\eqref{prob:cc general} with a joint chance constraint corresponding to the safety set $\mathcal{S}(\bm{x}) = \{\bm{\xi} \in \mathbb{R}^K \mid \bm{a}^\top_m \bm{x} < \bm{b}^\top_m\bm{\xi} + b_{m0} ~\forall m \in [M]\}$. As discussed in Section~\ref{sec:introduction}, we can assume that $\bm{b}_m \ne \bm{0}$ for all $m \in [M]$. By Proposition~\ref{prop:joint cc}, the distributionally robust chance constrained program~\eqref{prob:cc general} is thus equivalent to the deterministic optimization problem
\begin{equation}
\label{prob:joint cc reformulation abstract}
Z^\star_{\rm JCC}= \min_{(\bm{x},\bm{s}, t)\in \mathcal{C}_{\rm JCC}} \bm c^\top \bm x,
\end{equation}
whose feasible region is given by
\begin{equation*}
\mathcal{C}_{\rm JCC} = \left\{(\bm{x},\bm{s}, t) \in \mathcal{X} \times \mathbb{R}^N_+ \times \mathbb{R} ~\left|~
\begin{array}{ll}
\varepsilon N t - \mathbf{e}^\top\bm{s} \geq \theta N \\
\bigg(\displaystyle \min_{m \in [M]}
\bigg\{\dfrac{\bm{b}^\top_m\bmh{\xi}_i + b_{m0} - \bm{a}^\top_m\bm{x}}{\|\bm{b}_m\|_*} \bigg\}\bigg)^+ \geq t - s_i &~\forall i \in [N] \\
\end{array}
\right.\right\}.
\end{equation*}	
In analogy to Section~\ref{sec:cvar_individual_cc}, one can again show that any convex inner approximation of $\mathcal{C}_{\rm JCC}$ is weakly dominated by a polyhedron of the form
\begin{equation*}
\mathcal{C}_{\rm JCC}(\bm{\kappa}) = \left\{(\bm{x},\bm{s}, t) \in \mathcal{X} \times \mathbb{R}^N_+ \times \mathbb{R} ~\left|~
\begin{array}{ll}
\varepsilon N t - \mathbf{e}^\top\bm{s} \geq \theta N \\
\kappa_i \left(\displaystyle
\dfrac{\bm{b}^\top_m\bmh{\xi}_i + b_{m0} - \bm{a}^\top_m\bm{x}}{\|\bm{b}_m\|_*} \right) \geq t - s_i &~\forall m \in [M], ~i \in [N]  \\
\end{array}
\right.\right\}
\end{equation*}	
for some vector of slope parameters $\bm{\kappa} \in [0,1]^N$. The following assertion is akin to \Cref{prop:dominant convex approximation_1} and formalizes this statement. Its proof is omitted for the sake of brevity.

\begin{proposition}\label{prop:dominant convex approximation_JCC}
For any convex set $\mathcal{W} \subseteq \mathcal{C}_{\rm JCC}$, there exists $\bm{\kappa} \in[0,1]^N$ with $\mathcal{W} \subseteq \mathcal{C}_{\rm JCC} (\bm{\kappa}) \subseteq \mathcal{C}_{\rm JCC}$.
\end{proposition}

\Cref{prop:dominant convex approximation_JCC} implies that amongst all convex conservative approximations to problem~\eqref{prob:joint cc reformulation abstract} it is sufficient to consider the family of linear programs 
\begin{equation}
\label{prob:linear approximation joint cc_1}
Z^\star_{\rm JCC}(\bm{\kappa}) =\min_{(\bm{x},\bm{s}, t)\in \mathcal{C}_{\rm JCC}(\bm \kappa)} \bm c^\top \bm x
\end{equation}
parameterized by $\bm{\kappa} \in [0,1]^N$. One can show that the best approximation within this family is exact. The proof of this result is similar to that of Proposition~\ref{prop:exact dominant convex approximation_1} and thus omitted.
\begin{proposition}\label{prop:exact dominant convex approximation_2}
We have $Z^\star_{\rm JCC} = \min_{\bm\kappa \in [0,1]^N} Z^\star_{\rm JCC}(\bm{\kappa})$.
\end{proposition}

Unfortunately, finding the best slope parameters $\bm \kappa^\star\in [0,1]^N$ is again NP-hard, but optimizing over the subclass of uniform slope parameters $\bm{\kappa} = \kappa \mathbf{e}$ for  $\kappa \in [0,1]$ is easy, and $\bm{\kappa} = \mathbf{e}$ is optimal. This result is reminiscent of Proposition~\ref{prop:cvar best_1}, and thus its proof is omitted for the sake of brevity.

\begin{proposition}\label{prop:cvar best_2}
We have $\min_{\kappa \in [0,1]} Z^\star_{\rm JCC}(\kappa\mathbf{e})= Z^\star_{\rm JCC}(\mathbf{e})$.
\end{proposition}

We now demonstrate that $\mathcal{C}_{\rm JCC}(\mathbf{e})$ can again be interpreted as the feasible set of a worst-case CVaR constraint. To see this, denote by $\Delta_{++}^M =\{\bm w\in(0,1)^M \mid \mathbf{e}^\top\bm w=1 \}$ the relative interior of the probability simplex in $\mathbb{R}^M$ and observe that for any vector of scaling factors $\bm w\in\Delta_{++}^M$ we have
\begin{eqnarray*}
\mathbb{P}[\bmt{\xi} \in \mathcal{S}(\bm{x})] \ge 1-\varepsilon ~& \iff &~ \mathbb{P}\Big[\max_{m \in [M]}\{w_m(\bm{a}^\top_m\bm{x} - \bm{b}^\top_m\bm{\xi} - b_{m0})\} \geq 0\Big]\le\varepsilon \\
~&\iff &~ \mathbb{P}\text{-VaR}_{\varepsilon}\Big(\max_{m \in [M]}\{w_m(\bm{a}^\top_m\bm{x} - \bm{b}^\top_m\bm{\xi} - b_{m0})\}\Big) \leq 0 \\
~&\Longleftarrow &~ \mathbb{P}\text{-CVaR}_{\varepsilon}\Big(\max_{m \in [M]}\{w_m(\bm{a}^\top_m\bm{x} - \bm{b}^\top_m\bm{\xi} - b_{m0})\}\Big) \leq 0,
\end{eqnarray*}
where the first equivalence follows from the definition of the safety set $\mathcal S(\bm x)$. We emphasize that the exact reformulations of the joint chance constraint in the first two lines of the above expression are unaffected by the particular choice of $\bm w$ (that is, for any $\bm w,\bm w'\in \Delta_{++}^M$, a decision $\bm x$ is feasible for $\bm w$ if and only if it is feasible for $\bm w'$), while the CVaR approximation changes with $\bm w$. Thus, the quality of the CVaR approximation can be tuned by varying $\bm w\in\Delta^M_{++}$ (see \citealt{ordoudis2021energy} for an application of this tuning in energy and reserve dispatch). Note also that the overall normalization $\mathbf{e}^\top \bm w=1$ is non-restrictive because the CVaR is positive homogeneous.

We now introduce a family of worst-case CVaR constrained programs
\begin{equation}
\label{eq:wc-cvar-program_joint}
Z^\star_{\rm CVaR}(\bm{w}) = \left\{ 
\begin{array}{cl}
\displaystyle \min_{\bm{x} \in \mathcal{X}} & \bm{c}^\top\bm{x} \\
{\rm s.t.} & \displaystyle \sup_{\mathbb{P} \in \mathcal{F}(\theta)} \mathbb{P}\text{-CVaR}_{\varepsilon}\Big(\max_{m \in [M]}\{w_m(\bm{a}^\top_m\bm{x} - \bm{b}^\top_m\bm{\xi} - b_{m0})\}\Big) \leq 0
\end{array}\right.
\end{equation}
parameterized by $\bm w\in\Delta_{++}^M$, all of which conservatively approximate the ambiguous chance constrained program~\eqref{prob:cc general}, that is, $Z^\star_{\rm JCC} \le Z^\star_{\rm CVaR}(\bm{w})$. In fact, the family \eqref{eq:wc-cvar-program_joint} contains an instance that is equivalent to the best bounding problem of the form \eqref{prob:linear approximation joint cc_1} with uniform slope parameters.

\begin{proposition}\label{prop:joint cc worst-case CVaR approximation_1}
We have $Z^\star_{\rm CVaR}(\bm{w}^\star) = Z^\star_{\rm JCC}(\mathbf{e})$ for $\bm{w}^\star\in\Delta_{++}^M$ defined through
\[	
w^\star_m = \frac{\|\bm{b}_m\|_*^{-1}}{\sum_{\ell\in [M]} \|\bm b_{\ell}\|_*^{-1}} \quad \forall m \in [M].
\]
\end{proposition}

As the quality of the CVaR approximation in~\eqref{eq:wc-cvar-program_joint} depends on the choice of $\bm w$, it would be desirable to identify the best (least conservative)  approximation by solving $\min_{\bm w\in \Delta_{++}^M}Z^\star_{\rm CVaR}(\bm{w})$. This could be achieved, for instance, by treating $\bm w\in \Delta_{++}^M$ as an additional decision variable in~\eqref{eq:wc-cvar-program_joint}. Unfortunately, the resulting optimization problem involves bilinear terms in $\bm x$ and $\bm w$ and is thus non-convex. Finding the best CVaR approximation therefore appears to be computationally challenging. Even if the optimal scaling parameters were known, we will see in Section~\ref{sec:bonferroni} that the corresponding instance of problem~\eqref{eq:wc-cvar-program_joint} would generically provide a {\em strict} upper bound on $Z^\star_{\rm JCC}$.

The CVaR approximation~\eqref{eq:wc-cvar-program_joint} can again be interpreted as imposing a lower bound on the costs of moving training samples to the unsafe set. To see this, we define the \emph{minimum signed distance} between a point $\bm{\xi}$ and a family of closed sets $\mathcal{C}_m$, $m \in [M]$, as $\mathbf{min\textbf{-}dist} (\bm{\xi}, \{\mathcal{C}_m \}_{m \in [M]}) = \min_{m \in [M]} \, \mathbf{sgn\textbf{-}dist} (\bm{\xi}, \mathcal{C}_m)$. We then obtain the following result (\emph{cf.}~Proposition~\ref{prop:CVaR equivalence_1}).

\begin{proposition}\label{prop:min-signed-distances}
If $\bm w$ is set to $\bm w^\star$ as defined in Proposition~\ref{prop:joint cc worst-case CVaR approximation_1}, then we have
$$
\mspace{-5mu}
\sup_{\mathbb{P} \in \mathcal{F}(\theta)} \mathbb{P}\text{\emph{-CVaR}}_{\varepsilon}\Big(\max_{m \in [M]}\{w_m^\star (\bm{a}^\top_m\bm{x} - \bm{b}^\top_m\bm{\xi} - b_{m0})\}\Big) \leq 0
~\Longleftrightarrow~ 
\dfrac{1}{N}\sum_{i = 1}^{\varepsilon N} \mathbf{min\textbf{-}dist} (\bmh{\xi}_{\pi_i(\bm{x})}, \{ \mathcal{H}_m(\bm{x}) \}_{m \in [M]}) \ge \theta,
$$
where $\bm{\pi}(\bm{x})$ orders the data points $\bmh{\xi}_i$ by their minimum signed distances to the family of closed sets $\mathcal{H}_m(\bm{x}) = \{\bm{\xi} \in \mathbb{R}^K \mid \bm{a}^\top_m\bm{x} \geq \bm{b}^\top_m\bm{\xi} + b_{m0}\}$, $m \in [M]$.
\end{proposition}

The proof of \Cref{prop:min-signed-distances} closely resembles that of \Cref{prop:CVaR equivalence_1} and is therefore omitted.

\begin{corollary}\label{prop:joint cc CVaR is exact_1}
If $\bm w$ is set to $\bm w^\star$ as defined in \Cref{prop:joint cc worst-case CVaR approximation_1}, then the worst-case CVaR approximation is exact, that is, $Z^\star_{\rm CVaR}(\bm{w}^\star) = Z^\star_{\rm JCC}$, under either of the following conditions. 
\begin{itemize}
\item[(i)] We have $\bmh{\xi}_i \in \mathcal{S}(\bm{x}^\star)$ for all $i \in [N]$, where $\bm{x}^\star$ is optimal in~\eqref{prob:cc general}.
\item[(ii)] We have $\varepsilon \leq 1/N$. 
\end{itemize}
\end{corollary}

The proof is similar to that of \Cref{prop:individual cc CVaR is exact_1} and is thus omitted. 

\begin{remark}
In analogy to the individual chance constrained programs, Propositions~\ref{prop:exact dominant convex approximation_2},~\ref{prop:cvar best_2} and~\ref{prop:joint cc worst-case CVaR approximation_1} establish the best-in-class performance of the worst-case CVaR approximation among the inner linearizations $\mathcal{C}_{\textrm{JCC}}$ with uniform slopes, and Corollary~\ref{prop:joint cc CVaR is exact_1} identifies sufficient conditions under which the worst-case CVaR approximation is exact.
\end{remark}

\section{Bonferroni Approximation}\label{sec:bonferroni}

We next investigate the Bonferroni approximation, which applies to joint chance constrained programs with right-hand side uncertainty, that is, a subclass of problem~\eqref{prob:cc general} where $\mathcal{S}(\bm{x}) = \{\bm{\xi} \in \mathbb{R}^K \mid \bm{a}_m^\top \bm{x} < \bm{b}_m^\top \bm{\xi} + b_{m0} ~\forall m \in [M] \}$ for $\bm{a}_m \in \mathbb{R}^L$, $\bm{b}_m \in \mathbb{R}^K$ and $b_{m0} \in \mathbb{R}$, $m \in [M]$. In this context, note that Bonferroni's inequality (which is also known as the union bound) implies that
$$
\mathbb{P}[\bmt{\xi} \notin \mathcal{S}(\bm{x})] = \mathbb{P}[\bm{a}^\top_1 \bm{x} \ge \bm{b}^\top_1 \bmt{\xi} + b_{10} \quad \text{or}\quad \cdots \quad \text{or} \quad \bm{a}^\top_M \bm{x} \ge \bm{b}^\top_M \bmt{\xi} + b_{M0}] \leq \sum_{m \in [M]}\mathbb{P}[\bm{a}^\top_m \bm{x} \ge \bm{b}^\top_m \bmt{\xi} + b_{m0}].
$$
Taking the supremum over all distributions in the Wasserstein ball then yields the bound
\begin{equation}
\label{eq:bonferroni1}
\sup_{\mathbb P\in\mathcal F(\theta)} \mathbb{P}[\bmt{\xi} \notin \mathcal{S}(\bm{x})] \leq \sup_{\mathbb P\in\mathcal F(\theta)} \sum_{m \in [M]}\mathbb{P}[\bm{a}^\top_m \bm{x} \ge \bm{b}^\top_m \bmt{\xi} + b_{m0}] \leq \sum_{m \in [M]} \sup_{\mathbb P\in\mathcal F(\theta)} \mathbb{P}[\bm{a}^\top_m \bm{x} \ge \bm{b}^\top_m\bmt{\xi} + b_{m0}].
\end{equation}
For any collection of risk thresholds $\varepsilon_m \geq 0$, $m \in [M]$, such that $\sum_{m \in [M]} \varepsilon_m = \varepsilon$, the family of {\em individual} chance constraints
\begin{equation}
\label{eq:bonferroni2}
\sup_{\mathbb{P} \in \mathcal{F}(\theta)} \mathbb{P}[\bm{a}^\top_m \bm{x} \ge \bm{b}^\top_m \bmt{\xi} + b_{m0}] \leq \varepsilon_m \quad\forall m \in [M]
\end{equation}
thus provides a conservative approximation for the original {\em joint} chance constraint in~\eqref{prob:cc general} because
$$
\sup_{\mathbb P\in\mathcal F(\theta)} \mathbb{P}[\bmt{\xi} \notin \mathcal{S}(\bm{x})] \leq \sum_{m \in [M]} \sup_{\mathbb P\in\mathcal F(\theta)} \mathbb{P}[\bm{a}^\top_m \bm{x} \ge \bm{b}^\top_m \bmt{\xi} + b_{m0}]\leq \sum_{m \in [M]} \varepsilon_m = \varepsilon,
$$
where the two inequalities follow from~\eqref{eq:bonferroni1} and~\eqref{eq:bonferroni2}, respectively. We thus refer to~\eqref{eq:bonferroni2} as the {\em Bonferroni approximation} of the original chance constraint in problem~\eqref{prob:cc general}. The Bonferroni approximation is attractive because the individual chance constraints in~\eqref{eq:bonferroni2} are equivalent to simple linear inequalities. \mbox{To see this, note that each individual chance constraint in~\eqref{eq:bonferroni2} can be written as}
\begin{eqnarray*}
\sup_{\mathbb{P} \in \mathcal{F}(\theta)} \mathbb{P}[\bm{a}^\top_m \bm{x} \ge \bm{b}^\top_m \bmt{\xi} + b_{m0}] \leq \varepsilon_m 
~&\iff &~ \sup_{\mathbb{P} \in \mathcal{F}(\theta)} \mathbb{P}\text{-VaR}_{\varepsilon_m}(\bm{a}^\top_m\bm{x} - \bm{b}^\top_m\bmt{\xi} - b_{m0}) \leq 0 \\
~&\iff &~ \sup_{\mathbb{P} \in \mathcal{F}(\theta)} \mathbb{P}\text{-VaR}_{\varepsilon_m}(- \bm{b}^\top_m\bmt{\xi}) + \bm{a}^\top_m\bm{x} - b_{m0} \leq 0 \\
~&\iff &~ \sup_{\mathbb{P} \in \mathcal{F}(\theta)} \mathbb{P}\text{-VaR}_{\varepsilon_m}(- \bm{b}^\top_m\bmt{\xi}) \leq b_{m0} - \bm{a}^\top_m\bm{x},
\end{eqnarray*}
where the second equivalence holds because the value-at-risk is translation invariant. The $m^{\rm th}$~individual chance constraint in~\eqref{eq:bonferroni2} thus simplifies to the linear inequality $\bm{a}^\top_m\bm{x} \leq b_{m0} - \eta_m$, where the constant $\eta_m= \sup_{\mathbb{P} \in \mathcal{F}(\theta)} \mathbb{P}\text{-VaR}_{\varepsilon_m}(- \bm{b}^\top_m\bmt{\xi})$
is independent of $\bm{x}$ and can thus be computed offline. Specifically, by using Corollary~5.3 of \citet{Esfahani_Kuhn_2017}, we can express~$\eta_m$ as the optimal value of a deterministic optimization problem, that is,
$$
\eta_m=\left\{ \begin{array}{cll}
\displaystyle \min_{\bm{\alpha}, \beta, \bm{w}, \eta} & \eta \\
{\rm s.t.} & \theta\beta + \dfrac{1}{N} \displaystyle \sum_{i \in [N]} \alpha_i  \leq \varepsilon_m \\
& \alpha_i \geq 1 - w_i(\eta + \bm{b}^\top_m\bmh{\xi}_i)  &~\forall i \in [N] \\
& \beta \ge w_i \|\bm{b}_m\|_* &~\forall i \in [N] \\
& \bm \alpha\ge \bm 0,~\bm{w} \geq \bm{0}. 
\end{array}\right.
$$
The product of $\eta$ and $w_i$ in the second constraint group renders this problem non-convex. As the problem reduces to a linear program for any fixed value of the scalar decision variable $\eta$, however, $\eta_m$ can be computed efficiently to any accuracy by a line search along $\eta$. In summary, under the Bonferroni approximation the chance constrained program~\eqref{prob:cc general} thus reduces to a highly tractable linear program. However, the quality of the approximation relies on the choice of the individual risk thresholds $\{\varepsilon_m\}_{m \in [M]}$. It is often recommended to set $\varepsilon_m = \varepsilon/M$ for all $m \in [M]$, but \cite{Chen_Sim_Sun_2007} have shown that this choice can be conservative when the safety conditions are positively correlated. Optimizing over all admissible choices of $\{\varepsilon_m\}_{m \in [M]}$ is impractical because $\eta_m$ generically displays a non-convex dependence on~$\varepsilon_m$. Moreover, we will see that the Bonferroni approximation can be conservative even if the risk thresholds $\{\varepsilon_m\}_{m \in [M]}$ are chosen optimally.

\section{ALSO-X Approximation}\label{sec:ALSO-X}

Finally, we assume that the safety set is represented as $\mathcal{S}(\bm{x}) = \{\bm{\xi} \in \mathbb{R}^K \mid g(\bm{x}, \bm{\xi}) < 0\}$ for a generic function $g : \mathbb{R}^L \times \mathbb{R}^K \mapsto \mathbb{R}$. In particular, we recover individual chance constraints through the choice $g(\bm{x}, \bm{\xi}) = (\bm{A}\bm{\xi} + \bm{a})^\top \bm{x} - \bm{b}^\top\bm{\xi} - b_0$ (\emph{cf.}~Section~\ref{sec:cvar_individual_cc}), whereas joint chance constraints with right-hand side uncertainty correspond to $g(\bm{x}, \bm{\xi}) = \max_{m \in [M]} \{ \bm{a}^\top_m \bm{x} - \bm{b}^\top_m\bm{\xi} - b_{m0} \}$ (\emph{cf.}~Section~\ref{sec:cvar_joint_cc}).

The ALSO-X approximation \citep{ahmed2017nonanticipative, jiang2022also} traces the mapping
\begin{equation*}
\delta
\quad \mapsto \quad
\bm{x}^\star(\delta) \in \argmin \left\{
\sup_{\mathbb{P} \in \mathcal{F}(\theta)} \mathbb{E}_{\mathbb{P}}[(g(\bm{x}, \bmt{\xi}))^+]
\; : \;
\bm{x} \in \mathcal{X}, \; \bm{c}^\top\bm{x} \leq \delta
\right\},
\end{equation*}
where an arbitrary minimizer $\bm{x}^\star(\delta)$ is chosen whenever the minimization problem for $\delta$ admits multiple optimal solutions, and it employs a bisection search to determine the smallest value of $\delta$ for which the associated minimizer $\bm{x}^\star(\delta)$ satisfies the distributionally robust chance constraint:
\begin{equation*}
\sup_{\mathbb{P} \in \mathcal{F}(\theta)} \mathbb{P}[g(\bm{x}^\star(\delta), \bmt{\xi}) \geq 0] \leq \varepsilon.
\end{equation*}
If an incumbent solution $\bm{x}^\star(\delta)$ satisfies the chance constraint, then ALSO-X decreases $\delta$, otherwise the value of $\delta$ is increased. The ALSO-X approximation thus attempts to find a minimum cost solution $\bm{x}^\star(\delta)$, $\delta > 0$, that satisfies the worst-case chance constraint in problem~\eqref{prob:cc general} by replacing the chance constraint with the worst-case expected constraint violation $\sup_{\mathbb{P} \in \mathcal{F}(\theta)} \mathbb{E}_{\mathbb{P}}[(g(\bm{x}, \bmt{\xi}))^+]$.

For both individual chance constraints and joint chance constraints with right-hand side uncertainty, determining a solution $\bm{x}^\star(\delta)$ that minimizes the worst-case expected constraint violation can be reformulated as a deterministic and finite convex conic optimization problem; see, \textit{e.g.}, \citet[\S~5.1 and \S~7.1]{Esfahani_Kuhn_2017}. On the other hand, \Cref{thm:cc equivalent} allows us to efficiently verify whether $\bm{x}^\star(\delta)$ satisfies the distributionally robust chance constraint in problem~\eqref{prob:cc general}. The ALSO-X approximation thus requires the solution of a number of convex conic optimization problems and is therefore tractable.

\section{Incomparability of the Three Approximation Schemes}

We close this paper by demonstrating that for joint chance constrained programs with right-hand side uncertainty, the three approximation schemes are generally incomparable. To this end, we consider a data-driven setting where the number $N$ of available historical samples $\bmh{\xi}_1,\dots,\bmh{\xi}_N$ grows and where the radius $\theta$ of the Wasserstein ball decays to zero as $N$ approaches infinity. We also assume that in the limit, the ambiguity set $\mathcal{F}(\theta)$ contains the true data-generating distribution $\mathbb{P}_0$ with probability one. Both assumptions are standard in the literature, and they are satisfied by the usual choices of radii; see, \textit{e.g.}, \citet{blanchet2019robust} and \citet{kuhn2019wasserstein}. We also assume that $\|\cdot\|$ is a $p$-norm for some $p \in \mathbb{Q}$, $p > 1$.

\subsection{Incomparability between CVaR and Bonferroni Approximations}

We first provide two examples where either of the CVaR and the Bonferroni approximations is strictly less conservative than the other one. This is in stark contrast to Chebyshev ambiguity sets, where the worst-case CVaR approximation is known to dominate the Bonferroni approximation (see \citealt{Chen_Sim_Sun_Teo_2010} and \citealt{zymler2013distributionally}). 

\begin{example}\label{ex:data-driven_1}
Consider the following instance of the distributionally robust problem~\eqref{prob:cc general}:
\begin{equation}\label{data-driven_1}
\begin{array}{cll}
\displaystyle \min_{\bm{x}} & x_1 \\
{\rm s.t.} & \displaystyle \mathbb{P} [x_1 > \tilde{\xi}_1, ~x_2 > \tilde{\xi}_2] \geq 1-\varepsilon &~\forall \mathbb{P} \in \mathcal{F}(\theta) \\
& \underline{x}_1 \leq x_1 \leq \overline{x}_1, ~x_2 \geq 0.
\end{array}
\end{equation}
Here, we assume that $0 < \underline{x}_1 \leq \overline{x}_1 < 1$ and that the true data-generating distribution $\mathbb{P}_0$ is a two-point distribution which satisfies $\mathbb{P}_0 [(\tilde{\xi}_1, \tilde{\xi}_2) = (1,0)] = \rho$ and $\mathbb{P}_0 (\tilde{\xi}_1, \tilde{\xi}_2) = (0,0)] = 1-\rho$ for $\rho \in (0, 1)$.
\end{example}

\begin{proposition}\label{prop:bonferroni better than CVaR}
Let $\rho \in (\overline{x}_1 \varepsilon, \varepsilon)$. As $N \to \infty$, with probability going to $1$, we have that
\begin{enumerate}
\item[(i)] the Bonferroni approximation to~\eqref{data-driven_1} that replaces the joint chance constraint with
\begin{equation*}
\mathbb{P} [x_1 > \tilde{\xi}_1] \geq 1-\varepsilon_1 ~\forall \mathbb{P} \in \mathcal{F}(\theta), \quad
\mathbb{P} [x_2 > \tilde{\xi}_2] \geq 1-\varepsilon_2 ~\forall \mathbb{P} \in \mathcal{F}(\theta)
\end{equation*}
becomes exact if the risk thresholds $(\varepsilon_1, \varepsilon_2)$ are sufficiently close to $(\varepsilon, 0)$;
\item[(ii)] the worst-case CVaR approximation to~\eqref{data-driven_1} that replaces the joint chance constraint with
\begin{equation*}
\sup_{\mathbb{P} \in \mathcal{F}(\theta)}\mathbb{P}\text{\emph{-CVaR}}_{\varepsilon}\big(\max\big\{w_1(\tilde{\xi}_1 - x_1), w_2(\tilde{\xi}_2 - x_2)\big\}\big) \leq 0
\end{equation*}
becomes infeasible for any choice of scaling factors $(w_1, w_2) \in \Delta^2_{++}$.
\end{enumerate}
\end{proposition}

\begin{example}\label{ex:data-driven_2}
Consider the following instance of the distributionally robust problem~\eqref{prob:cc general}:
\begin{equation}\label{data-driven_2}
\begin{array}{cll}
\displaystyle \min_{\bm{x}} & x_3 \\
{\rm s.t.}&\mathbb{P}[x_1 > \tilde{\xi}, ~x_2 > \tilde{\xi}] \geq 1 - \varepsilon &~\forall \mathbb{P} \in \mathcal{F}(\theta) \\
&\underline{x} \leq x_1,x_2,x_3 \leq 1, ~x_3 \geq x_1, ~x_3 \geq x_2.
\end{array}
\end{equation}
Here, we assume $\frac{1}{2} < \underline{x} \leq 1$ and that the true data-generating distribution $\mathbb{P}_0$ is a two-point distribution which satisfies $\mathbb{P}_0[\tilde{\xi} = 1] = \rho$ and $\mathbb{P}_0[\tilde{\xi} = 0] = 1-\rho$ for $\rho \in (0,1)$.
\end{example}

\begin{proposition}\label{prop:CVaR better than bonferroni}
Let $\rho \in (\varepsilon/2, \underline{x}\varepsilon]$. As $N \to \infty$, with probability going to $1$, we have that
\begin{enumerate}
\item[(i)] the worst-case CVaR approximation to~\eqref{data-driven_2} that replaces the joint chance constraint with
\begin{equation*}
\sup_{\mathbb{P} \in \mathcal{F}(\theta)}\mathbb{P}\text{\emph{-CVaR}}_{\varepsilon}\big(\max\big\{w_1(\tilde{\xi} - x_1), w_2(\tilde{\xi} - x_2)\big\}\big) \leq 0
\end{equation*}
becomes exact if the scaling factors $(w_1, w_2)$ are $(\frac{1}{2}, \frac{1}{2})$; 
\item[(ii)] the Bonferroni approximation to~\eqref{data-driven_2} that replaces the joint chance constraint with
\begin{equation*}
\mathbb{P} [x_1 > \tilde{\xi}] \geq 1-\varepsilon_1 ~\forall \mathbb{P} \in \mathcal{F}(\theta), \quad
\mathbb{P} [x_2 > \tilde{\xi}] \geq 1-\varepsilon_2 ~\forall \mathbb{P} \in \mathcal{F}(\theta)
\end{equation*}
becomes infeasible for any choice of risk thresholds $(\varepsilon_1, \varepsilon_2)$.
\end{enumerate}	
\end{proposition}

\subsection{Incomparability between CVaR and ALSO-X Approximations}

We next provide two examples which demonstrate that the CVaR and the ALSO-X approximations are generically incomparable.

\begin{example}\label{example:CVaR_ALSO_1}\label{ex:data-driven_3}
Consider the following instance of the distributionally robust problem~\eqref{prob:cc general}:
\begin{equation}\label{data-driven_3}
\begin{array}{cll}
\displaystyle \min_{x} & x \\
{\rm s.t.} & \mathbb{P} [x > \tilde{\xi}_1, ~x < \tilde{\xi}_2] \geq 1-\varepsilon &~\forall \mathbb{P} \in \mathcal{F}(\theta) \\
&~ 0 \leq x \leq 1.
\end{array}
\end{equation}
Here, we assume that the true data-generating distribution $\mathbb{P}_0$ is a two-point distribution which satisfies $\mathbb{P}_0 [(\tilde{\xi}_1, \tilde{\xi}_2) = (1, -1)] = \rho$ and $\mathbb{P}_0 (\tilde{\xi}_1, \tilde{\xi}_2) = (-1, 1)] = 1-\rho$ for $\rho \in (0, 1)$.
\end{example}

\begin{proposition}\label{prop:ALSO-X better than CVaR}
Let $\rho \in (\varepsilon/2, \varepsilon)$. As $N \to \infty$, with probability going to $1$, we have that
\begin{enumerate}
\item[(i)] the ALSO-X approximation to~\eqref{data-driven_3} becomes exact;
\item[(ii)] the worst-case CVaR approximation to~\eqref{data-driven_3} that replaces the joint chance constraint with
\begin{equation*}
\sup_{\mathbb{P} \in \mathcal{F}(\theta)}\mathbb{P}\text{\emph{-CVaR}}_{\varepsilon}\big(\max\big\{w_1(\tilde{\xi}_1 - x), w_2(x - \tilde{\xi}_2)\big\}\big) \leq 0
\end{equation*}
becomes infeasible for any choice of scaling factors $(w_1, w_2) \in \Delta^2_{++}$.
\end{enumerate}	
\end{proposition}

\begin{example}\label{example:CVaR_ALSO_2}\label{ex:data-driven_4}
Consider the following instance of the distributionally robust problem~\eqref{prob:cc general}:
\begin{equation}\label{data-driven_4}
\begin{array}{cll}
\displaystyle \inf_{x} & -x \\
{\rm s.t.} & \displaystyle \mathbb{P} [x > \tilde{\xi}_1, ~x < \tilde{\xi}_2] \geq 1-\varepsilon &~\forall \mathbb{P} \in \mathcal{F}(\theta) \\
& 0 \leq x \leq 1. 
\end{array}
\end{equation}
Here, we assume that the true data-generating distribution $\mathbb{P}_0$ is a two-point distribution which satisfies $\mathbb{P}_0 (\tilde{\xi}_1, \tilde{\xi}_2) = (1,1)] = \rho$ and $\mathbb{P}_0 [(\tilde{\xi}_1, \tilde{\xi}_2) = (-1, 1)] = 1-\rho$ and  for $\rho \in (0, 1)$.
\end{example}

\begin{proposition}\label{prop:CVaR better than ALSO-X}
Let $\rho \in (0, \varepsilon/2)$. As $N \to \infty$, with probability going to $1$, we have that
\begin{enumerate}
\item[(i)] the worst-case CVaR approximation to~\eqref{data-driven_4} that replaces the joint chance constraint with
\begin{equation*}
\sup_{\mathbb{P} \in \mathcal{F}(\theta)}\mathbb{P}\text{\emph{-CVaR}}_{\varepsilon}\big(\max\big\{w_1(\tilde{\xi}_1 - x), w_2(x - \tilde{\xi}_2)\big\}\big) \leq 0
\end{equation*}
becomes exact if the scaling factors $(w_1, w_2)$ are $(\frac{1}{2}, \frac{1}{2})$; 
\item[(ii)] the ALSO-X approximation to~\eqref{data-driven_4} becomes infeasible.
\end{enumerate}	
\end{proposition}

\subsection{Incomparability between Bonferroni and ALSO-X Approximations}

Finally, we estabilish that either of the Bonferroni and the ALSO-X approximations can be strictly less conservative than the other one.

\begin{proposition}\label{prop:ALSO-X better than Bonferroni}
Let $\rho \in (\varepsilon/2, \varepsilon)$. As $N \to \infty$, with probability going to $1$, we have that the Bonferroni approximation to~\eqref{data-driven_3} that replaces the joint chance constraint with
\begin{equation*}
\mathbb{P} [x > \tilde{\xi}_1] \geq 1-\varepsilon_1 ~\forall \mathbb{P} \in \mathcal{F}(\theta), \quad
\mathbb{P} [x < \tilde{\xi}_2] \geq 1-\varepsilon_2 ~\forall \mathbb{P} \in \mathcal{F}(\theta)
\end{equation*}
becomes infeasible for any choice of risk thresholds $(\varepsilon_1, \varepsilon_2)$.
\end{proposition}

Recall from Proposition~\ref{prop:ALSO-X better than CVaR} that the ALSO-X approximation to~\eqref{data-driven_3} becomes exact under the setting of Proposition~\ref{prop:ALSO-X better than Bonferroni}.

\begin{proposition}\label{prop:bonferroni better than ALSO-X}
Let $\rho \in (0, \varepsilon/2)$. As $N \to \infty$, with probability going to $1$, we have that the Bonferroni approximation~\eqref{data-driven_4} that replaces the joint chance constraint with
\[
\mathbb{P} [x > \tilde{\xi}_1] \geq 1-\varepsilon_1 ~\forall \mathbb{P} \in \mathcal{F}(\theta), \quad
\mathbb{P} [x < \tilde{\xi}_2] \geq 1-\varepsilon_2 ~\forall \mathbb{P} \in \mathcal{F}(\theta)
\]
becomes exact if the risk thresholds $(\varepsilon_1, \varepsilon_2)$ are sufficiently close to $(\varepsilon, 0)$.
\end{proposition}

Recall from Proposition~\ref{prop:CVaR better than ALSO-X} that the ALSO-X approximation to~\eqref{data-driven_4} becomes infeasible under the setting of Proposition~\ref{prop:bonferroni better than ALSO-X}.

\section*{Acknowledgments}
The authors are grateful to Henry Lam and two anonymous referees for their thoughtful comments that substantially improved the paper. 
This research was supported by the ECS grant~CityU$\underline{~}$21502820, the SNSF grant BSCGI0$\underline{~}$157733 and the EPSRC grant EP/N020030/1. The authors thank Weijun Xie for helpful discussions on the ALSO-X approximation.

%%%%%%%%%%%%%%%%%%%%%%%%%%%%%%%% bibliography %%%%%%%%%%%%%%%%%%%%%%%%%%%%%%%%%%
\bibliographystyle{ormsv080}
\bibliography{References}

\begin{thebibliography}{22}
\expandafter\ifx\csname natexlab\endcsname\relax\def\natexlab#1{#1}\fi
\expandafter\ifx\csname url\endcsname\relax
  \def\url#1{{\tt #1}}\fi
\expandafter\ifx\csname urlprefix\endcsname\relax\def\urlprefix{URL }\fi
\expandafter\ifx\csname urlstyle\endcsname\relax
  \expandafter\ifx\csname doi\endcsname\relax
  \def\doi#1{doi:\discretionary{}{}{}#1}\fi \else
  \expandafter\ifx\csname doi\endcsname\relax
  \def\doi{doi:\discretionary{}{}{}\begingroup \urlstyle{rm}\Url}\fi \fi

\bibitem[{Ahmed et~al.(2017)Ahmed, Luedtke, Song, and
  Xie}]{ahmed2017nonanticipative}
Ahmed, Shabbir, James Luedtke, Yongjia Song, Weijun Xie. 2017.
\newblock Nonanticipative duality, relaxations, and formulations for
  chance-constrained stochastic programs.
\newblock {\it Mathematical Programming\/} {\bf 162}(1) 51--81.

\bibitem[{Ben-Tal and Nemirovski(2001)}]{Ben-tal_Nemirovski_book}
Ben-Tal, Aharon, Arkadi Nemirovski. 2001.
\newblock {\it Lectures on modern convex optimization: analysis, algorithms,
  and engineering applications\/}.
\newblock SIAM.

\bibitem[{Blanchet et~al.(2019)Blanchet, Kang, and Murthy}]{blanchet2019robust}
Blanchet, Jose, Yang Kang, Karthyek Murthy. 2019.
\newblock Robust {\uppercase{w}}asserstein profile inference and applications
  to machine learning.
\newblock {\it Journal of Applied Probability\/} {\bf 56}(3) 830--857.

\bibitem[{Chen et~al.(2010)Chen, Sim, Sun, and Teo}]{Chen_Sim_Sun_Teo_2010}
Chen, Wenqing, Melvyn Sim, Jie Sun, Chung-Piaw Teo. 2010.
\newblock From {\uppercase{cv}}a{\uppercase{r}} to uncertainty set:
  implications in joint chance-constrained optimization.
\newblock {\it Operations Research\/} {\bf 58}(2) 470--485.

\bibitem[{Chen et~al.(2007)Chen, Sim, and Sun}]{Chen_Sim_Sun_2007}
Chen, Xin, Melvyn Sim, Peng Sun. 2007.
\newblock A robust optimization perspective on stochastic programming.
\newblock {\it Operations Research\/} {\bf 55}(6) 1058--1071.

\bibitem[{Chen et~al.(2022)Chen, Kuhn, and
  Wiesemann}]{Chen_Kuhn_Wiesemann_2022}
Chen, Zhi, Daniel Kuhn, Wolfram Wiesemann. 2022.
\newblock Data-driven chance constrained programs over \uppercase{W}asserstein
  balls.  \emph{Forthcoming in Operations Research}.

\bibitem[{Hanasusanto et~al.(2015)Hanasusanto, Roitch, Kuhn, and
  Wiesemann}]{Hanasusanto_Roitch_Kuhn_Wiesemann_2015}
Hanasusanto, Grani~A, Vladimir Roitch, Daniel Kuhn, Wolfram Wiesemann. 2015.
\newblock A distributionally robust perspective on uncertainty quantification
  and chance constrained programming.
\newblock {\it Mathematical Programming\/} {\bf 151}(1) 35--62.

\bibitem[{Ho-Nguyen et~al.(2020)Ho-Nguyen, K{\i}l{\i}n{\c{c}}-Karzan,
  K{\"u}{\c{c}}{\"u}kyavuz, and Lee}]{ho2020strong}
Ho-Nguyen, Nam, Fatma K{\i}l{\i}n{\c{c}}-Karzan, Simge
  K{\"u}{\c{c}}{\"u}kyavuz, Dabeen Lee. 2020.
\newblock Strong formulations for distributionally robust chance-constrained
  programs with left-hand side uncertainty under \uppercase{W}asserstein
  ambiguity.  \emph{Forthcoming in INFORMS Journal on Optimization}.

\bibitem[{Ho-Nguyen et~al.(2022)Ho-Nguyen, K{\i}l{\i}n{\c{c}}-Karzan,
  K{\"u}{\c{c}}{\"u}kyavuz, and Lee}]{ho2021distributionally}
Ho-Nguyen, Nam, Fatma K{\i}l{\i}n{\c{c}}-Karzan, Simge
  K{\"u}{\c{c}}{\"u}kyavuz, Dabeen Lee. 2022.
\newblock Distributionally robust chance-constrained programs with right-hand
  side uncertainty under \uppercase{W}asserstein ambiguity.
\newblock {\it Mathematical Programming\/} {\bf 196} 641--672.

\bibitem[{Ho-Nguyen and Wright(2022)}]{ho2022adversarial}
Ho-Nguyen, Nam, Stephen Wright. 2022.
\newblock Adversarial classification via distributional robustness with
  \uppercase{W}asserstein ambiguity.  \emph{Forthcoming in Mathematical
  Programming}.

\bibitem[{Ji and Lejeune(2021)}]{ji2021data}
Ji, Ran, Miguel Lejeune. 2021.
\newblock Data-driven distributionally robust chance-constrained optimization
  with \uppercase{W}asserstein metric.
\newblock {\it Journal of Global Optimization\/} {\bf 79}(4) 779--811.

\bibitem[{Jiang and Xie(2022)}]{jiang2022also}
Jiang, Nan, Weijun Xie. 2022.
\newblock \uppercase{ALSO}-\uppercase{X} and \uppercase{ALSO}-\uppercase{X}+:
  better convex approximations for chance constrained programs.
  \emph{Forthcoming in Operations Research}.

\bibitem[{Kuhn et~al.(2019)Kuhn, Mohajerin~Esfahani, Nguyen, and
  Shafieezadeh-Abadeh}]{kuhn2019wasserstein}
Kuhn, Daniel, Peyman Mohajerin~Esfahani, Viet~Anh Nguyen, Soroosh
  Shafieezadeh-Abadeh. 2019.
\newblock Wasserstein distributionally robust optimization: theory and
  applications in machine learning.  \emph{Operations Research \& Management
  Science in the Age of Analytics (INFORMS)}, 130--166.

\bibitem[{Mohajerin~Esfahani and Kuhn(2018)}]{Esfahani_Kuhn_2017}
Mohajerin~Esfahani, Peyman, Daniel Kuhn. 2018.
\newblock Data-driven distributionally robust optimization using the
  {\uppercase{w}}asserstein metric: performance guarantees and tractable
  reformulations.
\newblock {\it Mathematical Programming\/} {\bf 171}(1-2) 1--52.

\bibitem[{Nemirovski(2012)}]{nemirovski2012safe}
Nemirovski, Arkadi. 2012.
\newblock On safe tractable approximations of chance constraints.
\newblock {\it European Journal of Operational Research\/} {\bf 219}(3)
  707--718.

\bibitem[{Ordoudis et~al.(2021)Ordoudis, Nguyen, Kuhn, and
  Pinson}]{ordoudis2021energy}
Ordoudis, Christos, Viet~Anh Nguyen, Daniel Kuhn, Pierre Pinson. 2021.
\newblock Energy and reserve dispatch with distributionally robust joint chance
  constraints.
\newblock {\it Operations Research Letters\/} {\bf 49}(3) 291--299.

\bibitem[{Shen and Jiang(2022)}]{shen2022chance}
Shen, Haoming, Ruiwei Jiang. 2022.
\newblock Chance-constrained set covering with \uppercase{W}asserstein
  ambiguity.  \textit{Forthcoming in Mathematical Programming}.

\bibitem[{Xie(2021)}]{xie2021distributionally}
Xie, Weijun. 2021.
\newblock On distributionally robust chance constrained programs with
  {\uppercase{w}}asserstein distance.
\newblock {\it Mathematical Programming\/} {\bf 186}(1) 115--155.

\bibitem[{Xie and Ahmed(2020)}]{xie2020bicriteria}
Xie, Weijun, Shabbir Ahmed. 2020.
\newblock Bicriteria approximation of chance-constrained covering problems.
\newblock {\it Operations Research\/} {\bf 68}(2) 516--533.

\bibitem[{Xie et~al.(2019)Xie, Ahmed, and Jiang}]{xie2019optimized}
Xie, Weijun, Shabbir Ahmed, Ruiwei Jiang. 2019.
\newblock Optimized {\uppercase{b}}onferroni approximations of distributionally
  robust joint chance constraints.
\newblock {\it Mathematical Programming\/} {\bf 191}(1) 79--112.

\bibitem[{Yang and Xu(2016)}]{yang2016distributionally}
Yang, Wenzhuo, Huan Xu. 2016.
\newblock Distributionally robust chance constraints for non-linear
  uncertainties.
\newblock {\it Mathematical Programming\/} {\bf 155}(1-2) 231--265.

\bibitem[{Zymler et~al.(2013)Zymler, Kuhn, and
  Rustem}]{zymler2013distributionally}
Zymler, Steve, Daniel Kuhn, Ber{\c{c}} Rustem. 2013.
\newblock Distributionally robust joint chance constraints with second-order
  moment information.
\newblock {\it Mathematical Programming\/} {\bf 137}(1-2) 167--198.

\end{thebibliography}
%%%%%%%%%%%%%%%%%%%%%%%%%%%%%%%% bibliography %%%%%%%%%%%%%%%%%%%%%%%%%%%%%%%%%%
%\begin{thebibliography}{}
%\end{thebibliography}
%%%%%%%%%%%%%%%%%
%%%%%%%%%%%%%%%%%

%%%%%%%%%%%%%%%%%%%%%%%%%%%%%%%% Appendix %%%%%%%%%%%%%%%%%%%%%%%%%%%%%%%%%%%%%%
\begin{appendices}
\section*{Proofs}

\noindent \emph{Proof of \Cref{prop:dominant convex approximation_1}.} $\;$
It is clear that $\mathcal{C}_{\rm ICC} (\bm{\kappa}) \subseteq \mathcal{C}_{\rm ICC}$ for every $\bm{\kappa} \in[0,1]^N$. Next, we show that for every $i \in [N]$ there exists $\kappa_i \in [0, 1]$ such that the constraint
$\kappa_i ((\bm{b} - \bm{A}^\top\bm{x})^\top\bmh{\xi}_i + b_0 - \bm{a}^\top\bm{x}) \geq t - s_i$ is valid for $\mathcal{W}$. The resulting set $\mathcal{C}_{\rm ICC}(\bm{\kappa})$ is thus a convex outer approximation of $\mathcal{W}$.

To determine $\kappa_i$, consider the sets
$\mathcal{W}_i = \{(\bm{x}, s_i, t) \mid (\bm{x}, \bm{s}, t) \in \mathcal{W}\}$
and $\mathcal{V}_i = \{(\bm{x}, s_i, t) \mid ((\bm{b} - \bm{A}^\top\bm{x})^\top\bmh{\xi}_i + b_0 - \bm{a}^\top\bm{x})^+ < t - s_i\}$.
By construction, $\mathcal{W}_i$ and $\mathcal{V}_i$ are intersection-free and convex. Thus, they admit a separating hyperplane. The same holds true if we replace $\mathcal{W}_i$ with
$$
\overline{\mathcal{W}}_i = \text{conv} \big(\mathcal{W}_i \cup \{(\bm{x}, s_i, t) \mid ((\bm{b} - \bm{A}^\top\bm{x})^\top\bmh{\xi}_i + b_0 - \bm{a}^\top\bm{x}, t - s_i) = (0, 0)\}\big).
$$
The separating hyperplane between $\overline{\mathcal{W}}_i$ and $\mathcal{V}_i$ must satisfy $t - s_i = 0$ whenever $(\bm{b} - \bm{A}^\top\bm{x})^\top\bmh{\xi}_i + b_0 - \bm{a}^\top\bm{x} = 0$. In other words, the separating hyperplane must be of the form $\kappa_i((\bm{b} - \bm{A}^\top\bm{x})^\top\bmh{\xi}_i + b_0 - \bm{a}^\top\bm{x}) = t - s_i$ for some $\kappa_i\in [0,1]$. Thus, the claim follows.
\hfill \Halmos
\vspace{2mm}

\noindent \emph{Proof of \Cref{prop:exact dominant convex approximation_1}.} $\;$
It follows from \Cref{prop:individual cc} that
\begin{equation}
\label{eq:aux1}
\min_{\bm \kappa\in [0,1]^N} Z^\star_{\rm ICC}(\bm{\kappa})=\left\{ \begin{array}{cll}
\displaystyle \min_{\bm{s}, t, \bm{x},\bm \kappa} & \bm{c}^\top\bm{x} \\
{\rm s.t.} & \varepsilon N t - \mathbf{e}^\top\bm{s} \geq \theta N  \|\bm{b} - \bm{A}^\top\bm{x}\|_*\\
& \kappa_i((\bm{b} - \bm{A}^\top\bm{x})^\top\bmh{\xi}_i + b_0 - \bm{a}^\top\bm{x}) \geq t - s_i &~\forall i \in [N] \\
& \bm{s} \geq \bm{0}, ~\bm{x} \in \mathcal{X},~ \bm \kappa\in [0,1]^N.
\end{array}\right.
\end{equation}
For any fixed $(\bm{x}, \bm{s}, t)$, the optimal (that is, least restrictive) choice of $\bm \kappa$ satisfies
\begin{equation}
\label{eq:aux2}
\kappa_i = \left\{
\begin{array}{ll}
1 &\text{~if~} (\bm{b} - \bm{A}^\top\bm{x})^\top\bmh{\xi}_i + b_0 - \bm{a}^\top\bm{x} \geq 0, \\
0 &\text{~otherwise} 
\end{array}
\right.
\quad \forall i\in [N].
\end{equation}
Eliminating $\bm \kappa$ from~\eqref{eq:aux1} by substituting~\eqref{eq:aux2} into~\eqref{eq:aux1} converts the second constraint group to
\[
((\bm{b} - \bm{A}^\top\bm{x})^\top\bmh{\xi}_i + b_0 - \bm{a}^\top\bm{x})^+ \geq t - s_i \quad \forall i \in [N],
\]
which shows that~\eqref{eq:aux1} is equivalent to~\eqref{prob:individual cc reformulation abstract}. Thus, the claim follows.
\hfill \Halmos
\vspace{2mm}

\noindent \emph{Proof of \Cref{prop:cvar best_1}.} $\;$ We first show that problem~\eqref{prob:linear approximation individual cc_1} is infeasible for $\bm \kappa=\bm 0$, that is, $Z^\star_{\rm ICC}(\bm{0})=\infty$. Indeed, by the definition of $\mathcal{C}_{\rm ICC}(\bm{0})$ we have
\[
Z^\star_{\rm ICC}(\bm{0})=\left\{\begin{array}{cll}
\displaystyle \min_{\bm{s},t,\bm{x}} & \bm{c}^\top\bm{x} \\
{\rm s.t.} & \varepsilon N t - \mathbf{e}^\top\bm{s} \geq  \theta N \|\bm{b} - \bm{A}^\top\bm{x}\|_* \\
&\bm s\ge t\mathbf{e}, ~ \bm{s} \geq \bm{0},~\bm{x} \in \mathcal{X}.
\end{array}\right.
\]
Any feasible solution of the above problem satisfies $\varepsilon N t - \mathbf{e}^\top\bm{s}\le \varepsilon N t - N\max\{t,0\}\le 0$, where the first inequality follows from the constraints $\bm s\ge t\mathbf{e}$ and $ \bm{s} \geq \bm{0}$. As $\theta>0$, the constraint $\varepsilon N t - \mathbf{e}^\top\bm{s} \geq  \theta N \|\bm{b} - \bm{A}^\top\bm{x}\|_*$ is thus satisfied only if $\bm s=\bm 0$, $t=0$ and $\bm{A}^\top\bm{x}=\bm b$. However, the last equality contradicts our standing assumption that  $\bm{A}^\top\bm{x} \ne \bm{b}$ for all $\bm{x} \in \mathcal{X}$, confirming that the above problem is infeasible and~$Z^\star_{\rm ICC}(\bm{0})=\infty$. Thus, $Z^\star_{\rm ICC}(\kappa\mathbf{e})$ is minimized by some $\kappa \in (0,1]$.

If $\kappa\in (0,1]$, we can use the variable substitution $t\leftarrow \kappa t$ and $\bm{s} \leftarrow \kappa \bm{s}$ to re-express problem~\eqref{prob:linear approximation individual cc_1}~as
\[
Z^\star_{\rm ICC}(\kappa\mathbf{e}) =\left\{\begin{array}{cll}
\displaystyle \min_{\bm{s},t,\bm{x}} & \bm{c}^\top\bm{x} \\
{\rm s.t.} & \varepsilon N t - \mathbf{e}^\top\bm{s} \geq \dfrac{\theta N}{\kappa} \|\bm{b} - \bm{A}^\top\bm{x}\|_* \\
& (\bm{b} - \bm{A}^\top\bm{x})^\top\bmh{\xi}_i + b_0 - \bm{a}^\top\bm{x} \geq t - s_i &~\forall i \in [N] \\
& \bm{s} \geq \bm{0}, ~\bm{x} \in \mathcal{X}.
\end{array}\right.
\]
From this formulation it is evident that $\kappa^\star = 1$ is the best (least restrictive) choice of $\kappa \in (0,1]$. 
\hfill \Halmos
\vspace{2mm}

\noindent \emph{Proof of \Cref{prop:individual cc worst-case CVaR approximation_1}.} $\;$
Using now standard techniques, the worst-case CVaR in~\eqref{eq:wc-cvar-program} can be re-expressed as the optimal value of a finite conic program,
$$
\sup_{\mathbb{P} \in \mathcal{F}(\theta)} \mathbb{P}\text{-CVaR}_{\varepsilon}(\bm{a}^\top\bm{x} - b_0 + (\bm{A}^\top\bm{x}-\bm{b})^\top\bmt{\xi}) = \left\{\begin{array}{cll}
\displaystyle{\min_{\bm \alpha, \beta,\tau}} & \displaystyle \tau + \dfrac{1}{\varepsilon} \Big(\theta \beta + \dfrac{1}{N} \sum_{i \in [N]} \alpha_i\Big)  \\
{\rm s.t.} & \alpha_i  \geq \bm{a}^\top\bm{x} - b_0 + (\bm{A}^\top\bm{x}-\bm{b})^\top\bmh{\xi}_i -\tau &~\forall i \in [N] \\
& \bm{\alpha} \geq \bm{0}, ~ \beta \geq \|\bm{A}^\top\bm{x} - \bm{b}\|_*, 
\end{array}\right.
$$
see \citet[\S~5.1 and \S~7.1]{Esfahani_Kuhn_2017} for a detailed derivation. Substituting this reformulation into the worst-case CVaR constrained program~\eqref{eq:wc-cvar-program} yields
\[
Z^\star_{\rm CVaR}=\left\{ \begin{array}{cl}
\displaystyle \min_{\bm{x},\bm \alpha, \beta,\tau} & \bm{c}^\top\bm{x} \\
{\rm s.t.} & \tau + \dfrac{1}{\varepsilon} \Big(\theta \beta + \dfrac{1}{N} \displaystyle\sum_{i \in [N]} \alpha_i\Big) \leq 0\\
& \alpha_i  \geq \bm{a}^\top\bm{x} - b_0 + (\bm{A}^\top\bm{x}-\bm{b})^\top\bmh{\xi}_i -\tau \quad\forall i \in [N] \\
& \bm{\alpha} \geq \bm{0}, ~ \beta \geq \|\bm{A}^\top\bm{x} - \bm{b}\|_*,~ \bm{x} \in \mathcal{X}.
\end{array}\right.
\]
As $\theta > 0$ and $\varepsilon>0$, it is clear that $\beta = \|\bm{A}^\top\bm{x}-\bm{b}\|_*$ at optimality, and this insight allows us to eliminate~$\beta$ from the above optimization problem. Multiplying the first constraint with the positive constant $\varepsilon N$ while renaming $\bm \alpha$ as $\bm s$ and $\tau$ as $-t$ then shows that $Z^\star_{\rm CVaR}= Z^\star_{\rm ICC}(\mathbf{e})$.
\hfill \Halmos
\vspace{2mm}

\noindent \emph{Proof of \Cref{prop:CVaR equivalence_1}.} $\;$
It follows from the proof of \Cref{prop:individual cc worst-case CVaR approximation_1} that the worst-case CVaR constraint $\sup_{\mathbb{P} \in \mathcal{F}(\theta)}\mathbb{P}\text{-CVaR}_{\varepsilon}(\bm{a}^\top\bm{x} - b_0 + (\bm{A}^\top\bm{x}-\bm{b})^\top\bmt{\xi}) \leq 0$ holds if and only if
$$
0 \geq \left\{\begin{array}{cll}
\displaystyle{\min_{\bm{\alpha}, \tau}} & \displaystyle \tau + \dfrac{1}{\varepsilon} \Big(\theta \|\bm{A}^\top\bm{x} - \bm{b}\|_* + \dfrac{1}{N} \sum_{i \in [N]} \alpha_i\Big)  \\
{\rm s.t.} & \alpha_i  \geq \bm{a}^\top\bm{x} - b_0 + (\bm{A}^\top\bm{x}-\bm{b})^\top\bmh{\xi}_i -\tau &~\forall i \in [N] \\
& \bm{\alpha} \geq \bm{0},
\end{array}\right.
$$
which, by multiplying the objective function by the positive constant $\varepsilon N$ while renaming $\bm \alpha$ as $\bm s$ and $\tau$ as $-t$, is equivalent to
\begin{equation}\label{eq:cvar_ref_1}
\exists \bm{s} \geq \bm{0}, \, t \in \mathbb{R} \, : \,
\left\{
\begin{array}{ll}
\varepsilon N t - \mathbf{e}^\top\bm{s} \geq \theta N \|\bm{b} - \bm{A}^\top\bm{x}\|_* \\
(\bm{b} - \bm{A}^\top\bm{x})^\top\bmh{\xi}_{\pi_i(\bm{x})} + b_0 - \bm{a}^\top\bm{x} \geq t - s_i &~\forall i \in [N].
\end{array}
\right.
\end{equation}
This constraint system is satisfiable by $t \in \mathbb{R}$ and \emph{some} $\bm{s} \geq \bm{0}$ if and only if it is satisfiable by $t$ and $\bm{s}^\star (t)$ defined by $s_i^\star (t) = (t - ((\bm{b} - \bm{A}^\top\bm{x})^\top\bmh{\xi}_{\pi_i(\bm{x})} + b_0 - \bm{a}^\top\bm{x}))^+$, $i \in [N]$. Since the second constraint in~\eqref{eq:cvar_ref_1} is automatically satisfied by $\bm{s}^\star (t)$, we thus conclude that~\eqref{eq:cvar_ref_1} holds if and only if
\begin{align}
& \exists t \in \mathbb{R} \, : \, \varepsilon N t - \sum_{i \in [N]} (t - ((\bm{b} - \bm{A}^\top\bm{x})^\top\bmh{\xi}_{\pi_i(\bm{x})} + b_0 - \bm{a}^\top\bm{x}))^+ \geq \theta N \|\bm{b} - \bm{A}^\top\bm{x}\|_* \nonumber \\
\Longleftrightarrow \;\; & \max_{t \in \mathbb{R}} \bigg\{\varepsilon N t - \sum_{i \in [N]} (t - ((\bm{b} - \bm{A}^\top\bm{x})^\top\bmh{\xi}_{\pi_i(\bm{x})} + b_0 - \bm{a}^\top\bm{x}))^+ \bigg\} \geq \theta N \|\bm{b} - \bm{A}^\top\bm{x}\|_*. \label{eq:cvar_ref_2}
\end{align} 
The objective function of the embedded maximization problem on the left-hand side of~\eqref{eq:cvar_ref_2} is piecewise affine and concave in $t$. Moreover, by construction of $\bm{\pi} (\bm{x})$, we have
\begin{equation*}
(\bm{b} - \bm{A}^\top\bm{x})^\top\bmh{\xi}_{\pi_i(\bm{x})} + b_0 - \bm{a}^\top\bm{x}
\;\; \leq \;\;
(\bm{b} - \bm{A}^\top\bm{x})^\top\bmh{\xi}_{\pi_j(\bm{x})} + b_0 - \bm{a}^\top\bm{x}
\qquad \forall 1 \leq i \leq j \leq N.
\end{equation*}
The first-order optimality condition for non-smooth optimization then implies that the maximum on the left-hand side of~\eqref{eq:cvar_ref_2} is attained by $t^\star = (\bm{b} - \bm{A}^\top\bm{x})^\top\bmh{\xi}_{\pi_{\lfloor\varepsilon N\rfloor + 1}(\bm{x})} + b_0 - \bm{a}^\top\bm{x}$, which results in the equivalent constraint
\begin{equation*}
\sum_{i = 1}^{\varepsilon N} ((\bm{b} - \bm{A}^\top\bm{x})^\top\bmh{\xi}_{\pi_i(\bm{x})} + b_0 - \bm{a}^\top\bm{x}) \ge \theta N \|\bm{b} - \bm{A}^\top\bm{x}\|_*.
\end{equation*}
The result now follows if we divide both sides of the constraint by $N \|\bm{b} - \bm{A}^\top\bm{x}\|_*$.
\hfill \Halmos
\vspace{2mm}

\noindent \emph{Proof of \Cref{prop:individual cc CVaR is exact_1}.} $\;$
The first condition immediately follows from \Cref{thm:cc equivalent} and Proposition~\ref{prop:CVaR equivalence_1} since $\mathbf{sgn\textbf{-}dist} (\bmh{\xi}_i, \bar{\mathcal{S}} (\bm{x})) = \mathbf{dist} (\bmh{\xi}_i, \bar{\mathcal{S}} (\bm{x}))$ whenever $\bmh{\xi}_i \in \mathcal{S} (\bm{x})$. The second condition guarantees that $\bmh{\xi}_i \in \mathcal{S} (\bm{x})$, $i \in [N]$, for
any solution $\bm{x} \in \mathcal{X}$ that satisfies the ambiguous chance constraint in problem~\eqref{prob:cc general}. This, in turn, implies that the first condition of the corollary is satisfied as well.
\hfill \Halmos
\vspace{2mm}

\noindent \emph{Proof of \Cref{prop:joint cc worst-case CVaR approximation_1}.} $\;$
Using techniques introduced by \cite{Esfahani_Kuhn_2017}, the worst-case CVaR in~\eqref{eq:wc-cvar-program_joint} can be re-expressed as
\begin{align*}
& \sup_{\mathbb{P} \in \mathcal{F}(\theta)} \mathbb{P}\text{-CVaR}_{\varepsilon}\Big(\max_{m \in [M]}\{w_m(\bm{a}^\top_m\bm{x} - \bm{b}^\top_m\bm{\xi} - b_{m0})\}\Big) \\
& \qquad = \left\{\begin{array}{cll}
\displaystyle{\min_{\bm \alpha, \beta,\tau}} & \displaystyle \tau + \dfrac{1}{\varepsilon} \Big(\theta \beta + \dfrac{1}{N} \sum_{i \in [N]} \alpha_i\Big)  \\
{\rm s.t.} & \alpha_i  \geq w_m(\bm{a}^\top_m\bm{x} - \bm{b}^\top_m\bmh{\xi}_i - b_{m0}) -\tau &~\forall m \in [M], ~i \in [N] \\
& \beta \geq w_m\|\bm{b}_m\|_* &~\forall m \in [M] \\
& \bm{\alpha} \geq \bm{0}.
\end{array}\right.
\end{align*}
Substituting this reformulation into~\eqref{eq:wc-cvar-program_joint} yields
\[
Z^\star_{\rm CVaR}(\bm{w}) =\left\{ 
\begin{array}{cll}
\displaystyle \min_{\bm{x},\bm \alpha, \beta,\tau} & \bm{c}^\top\bm{x} \\
{\rm s.t.} & \tau + \dfrac{1}{\varepsilon} \Big(\theta \beta + \dfrac{1}{N} \displaystyle\sum_{i \in [N]} \alpha_i\Big) \leq 0\\
& \alpha_i \geq w_m(\bm{a}^\top_m\bm{x} - \bm{b}^\top_m\bmh{\xi}_i - b_{m0}) -\tau &~\forall m \in [M], ~i \in [N]\\
& \beta \geq w_m\|\bm{b}_m\|_* &~\forall m \in [M] \\
& \bm{\alpha} \geq \bm{0}, ~\bm{x} \in \mathcal{X}.
\end{array}\right.
\]
As $\theta > 0$ and $\varepsilon>0$, it is clear that $\beta = \max_{m \in [M]} \{w_m\|\bm{b}_m\|_*\}$ at optimality, and this insight allows us to eliminate~$\beta$ from the above optimization problem. Multiplying the first constraint by the positive constant $\varepsilon N/\beta$ and the second constraint group by the positive constant $1/\beta$ while applying the variable substitutions $\bm{s} \leftarrow \bm{\alpha}/\beta$ and $-t \leftarrow \tau/\beta$, we obtain
\begin{equation}
\label{eq:cvar-aux}
Z^\star_{\rm CVaR}(\bm{w}) =\left\{ 
\begin{array}{cll}
\displaystyle \min_{\bm s, t, \bm{x}} & \bm{c}^\top\bm{x} \\
{\rm s.t.} & \varepsilon N t - \mathbf{e}^\top\bm{s} \geq \theta N \\
& \dfrac{w_m\|\bm{b}_m\|_*}{\max\limits_{m \in [M]} \{w_m\|\bm{b}_m\|_*\}}\dfrac{(\bm{b}^\top_m\bmh{\xi}_i + b_{m0} - \bm{a}^\top_m\bm{x})}{\|\bm{b}_m\|_*} \geq t - s_i &~\forall m \in [M], ~i \in [N] \\
& \bm{s} \geq \bm{0}, ~\bm{x} \in \mathcal{X}.
\end{array}\right.
\end{equation}
Replacing $\bm{w}$ with $\bm{w}^\star$, the second constraint group in problem~\eqref{eq:cvar-aux} simplifies to
$$
\min_{m \in [M]} \Bigg\{\dfrac{(\bm{b}^\top_m\bmh{\xi}_i + b_{m0} - \bm{a}^\top_m\bm{x})}{\|\bm{b}_m\|_*}\Bigg\} \geq t - s_i ~\forall i \in [N],
$$
which reveals that the feasible set of problem~\eqref{eq:cvar-aux} coincides with $C_{\rm JCC}(\mathbf{e})$. This observation implies the postulated assertion that $Z^\star_{\rm CVaR}(\bm{w}^\star) = Z^\star_{\rm JCC}(\mathbf{e})$.
\hfill \Halmos
\vspace{2mm}

In the proofs of Propositions~\ref{prop:bonferroni better than CVaR}--\ref{prop:bonferroni better than ALSO-X} we will denote by $I$ the number of historical samples that coincide with the first atom of the data-generating distribution $\mathbb{P}_0$, that is, with $(1, 0)$ in Example~\ref{ex:data-driven_1}, with $1$ in Example~\ref{ex:data-driven_2}, with $(1, -1)$ in Example~\ref{ex:data-driven_3} and with $(1, 1)$ in Example~\ref{ex:data-driven_4}, respectively. We emphasize that $I$ is a random variable that depends on the number $N$ of samples. To avoid clutter, however, we do not use a tilde sign and suppress the dependence on $N$.
\vspace{2mm}

\noindent \emph{Proof of \Cref{prop:bonferroni better than CVaR}.} $\;$ 
We proceed in three steps. We first derive the optimal value of the classical chance constrained program associated with~\eqref{data-driven_1} under the true data-generating distribution $\mathbb{P}_0$ (Step~1). This value serves as a lower bound on the optimal value of problem~\eqref{data-driven_1}. We then show that with probability going to $1$ (w.p.~1) as $N \to \infty$, the Bonferroni approximation achieves this bound~(Step~2), whereas the worst-case CVaR approximation becomes infeasible~(Step~3).

{\em Step~1.} Since $\rho < \varepsilon$, the feasible region of the classical chance constrained program
\begin{equation*}
\begin{array}{cl}
\displaystyle \min_{\bm{x}} & x_1 \\
{\rm s.t.} & \displaystyle \mathbb{P}_0 [x_1 > \tilde{\xi}_1, ~x_2 > \tilde{\xi}_2] \geq 1-\varepsilon \\
& \underline{x}_1 \leq x_1 \leq \overline{x}_1, ~x_2 \geq 0
\end{array}
\end{equation*}
under the true data-generating distribution $\mathbb{P}_0$ is $\{ (x_1, x_2) \in \mathbb{R}^2 \mid x_1 \in [\underline{x}_1, \overline{x}_1], ~x_2 > 0 \}$. Hence, the optimal value of this problem is $\underline{x}_1>0$, which is attained by any $(x_1, x_2) \in \{ \underline{x}_1 \} \times (\mathbb{R}_+ \setminus \{ 0 \})$.

{\em Step~2.} 
Fix any $(x_1,x_2) \in [\underline{x}_1, \overline{x}_1]\times (\mathbb R_+\backslash\{0\})$, and denote by $\mathcal{S}_1(\bm{x}) = \{\bm{\xi} \mid x_1 > \xi_1 \}$ and $\mathcal{S}_2(\bm{x}) = \{\bm{\xi} \mid x_2 > \xi_2 \}$ the two safety sets of the Bonferroni approximation. If $\bmh{\xi}_i = (1, 0)^\top$, then $\bmh{\xi}_i \in \bar{\mathcal{S}}_1(\bm{x}) \cap \mathcal{S}_2(\bm{x})$ with $\mathbf{dist}(\bmh{\xi}_i, \bar{\mathcal{S}}_1(\bm{x})) = 0$ and $\mathbf{dist}(\bmh{\xi}_i, \bar{\mathcal{S}}_2(\bm{x})) = x_2$. Likewise, if $\bmh{\xi}_i = (0, 0)^\top$, then $\bmh{\xi}_i \in \mathcal{S}_1(\bm{x}) \cap \mathcal{S}_2(\bm{x})$ with $\mathbf{dist}(\bmh{\xi}_i, \bar{\mathcal{S}}_1(\bm{x})) = x_1$ and $\mathbf{dist}(\bmh{\xi}_i, \bar{\mathcal{S}}_2(\bm{x})) = x_2$. Under the appropriate permutations $\bm{\pi}^1 (\bm{x})$ and $\bm{\pi}^2 (\bm{x})$, \Cref{thm:cc equivalent} then implies that $\bm{x}$ satisfies both chance constraints of the Bonferroni approximation if and only if
\begin{equation}\label{eq:bonferroni_times_two}
\mspace{-10mu}
\dfrac{1}{N} \sum_{i = 1}^{\varepsilon_1 N}\mathbf{dist}(\bmh{\xi}_{\pi^1_i(\bm{x})}, \bar{\mathcal{S}}_1(\bm{x})) = \dfrac{1}{N}(\varepsilon_1 N - I)^+ x_1 \geq \theta
\quad \text{and} \quad
\dfrac{1}{N} \sum_{i = 1}^{\varepsilon_2 N}\mathbf{dist}(\bmh{\xi}_{\pi^2_i(\bm{x})}, \bar{\mathcal{S}}_2(\bm{x})) = \varepsilon_2 x_2 \geq \theta,
\end{equation}
where $I$ denotes the number of samples $\bmh{\xi}_i$, $i \in [N]$, that satisfy $\bmh{\xi}_i = (1, 0)^\top$.

Choose $\varepsilon_1 \in (p, \varepsilon)$ and $\varepsilon_2 = \varepsilon - \varepsilon_1$, as well as $x_1 = \underline{x}_1$ and any $x_2 \geq \theta/\varepsilon_2$. This choice of $(\varepsilon_1, \varepsilon_2)$ and $\bm{x}$ satisfies the second constraint in~\eqref{eq:bonferroni_times_two} by construction. To see that the first constraint in~\eqref{eq:bonferroni_times_two} is also satisfied w.p.~1 as $N \to \infty$, we note that $\frac{1}{N}(\varepsilon_1 N - I)^+ x_1 = (\varepsilon_1 - I/N)^+ x_1 \to (\varepsilon_1 - \rho)^+ x_1$ w.p.~1 as $N \to \infty$ due to the strong law of large numbers. We thus conclude that $\frac{1}{N}(\varepsilon_1 N - I)^+ x_1 > 0$ w.p.~1 as $N \to \infty$, and thus this quantity will exceed $\theta$, which goes to zero as $N$ approaches infinity.

{\em Step~3.} 
Given decision $(x_1, x_2)$, the left-hand side of the worst-case CVaR approximation to~\eqref{data-driven_1},
\begin{equation*}
\sup_{\mathbb{P} \in \mathcal{F}(\theta)} \mathbb{P}\text{-CVaR}_{\varepsilon}\big(\max\{w_1(\tilde{\xi}_1 - x_1), \; w_2(\tilde{\xi}_2 - x_2)\}\big),
\end{equation*}
can be expressed as the optimal value of the optimization problem
\begin{equation*}
\begin{array}{cll}
\displaystyle{\min_{\bm \alpha, \beta,\tau}} & \displaystyle \tau + \dfrac{1}{\varepsilon} \Big(\theta \beta + \dfrac{1}{N} \sum_{i \in [N]} \alpha_i\Big)  \\
{\rm s.t.} & \alpha_i  \geq w_1(\hat{\xi}_{i,1} - x_1) -\tau, \;\; \alpha_i  \geq w_2(\hat{\xi}_{i,2} - x_2) -\tau &~\forall i \in [N] \\
& \beta \geq w_1\|(1,0)^\top\|_*, \;\; \beta \geq w_2\|(0,1)^\top\|_* \\
& \bm{\alpha} \geq \bm{0}.
\end{array}
\end{equation*}
As $\theta > 0$, $\varepsilon > 0$ and $N > 0$, this problem is minimized by $\alpha_i^\star = (\max\{w_1(\hat{\xi}_{i,1}-x_1), w_2(\hat{\xi}_{i,2} - x_2)\} - \tau)^+$ and $\beta^\star = \max\{w_1, w_2\}$. Therefore, the worst-case CVaR approximation is feasible if and only if 
\[
\exists \tau \in \mathbb{R}: \;
0 \geq  \tau + \dfrac{1}{\varepsilon} \Big(\theta \cdot \max\{w_1, w_2\} + \dfrac{1}{N} \sum_{i \in [N]} (\max\{w_1(\hat{\xi}_{i,1} - x_1), w_2(\hat{\xi}_{i,2} - x_2)\} - \tau)^+\Big).
\]
Multiplying both sides by $\frac{\varepsilon N}{\max\{w_1, w_2\}}$ and applying the variable substitution $-t \leftarrow \frac{\tau}{\max\{w_1, w_2\}}$, the above condition becomes
\begin{align}
& \displaystyle \exists t \in \mathbb{R}: \;
0 \geq -\varepsilon N t + \theta N + \sum_{i \in [N]} \Bigg(\max\Bigg\{\frac{w_1(\hat{\xi}_{i,1}-x_1)}{\max\{w_1, w_2\}}, \frac{w_2(\hat{\xi}_{i,2}-x_2)}{\max\{w_1, w_2\}}\Bigg\}+ t\Bigg)^+ \nonumber \\
\Longleftrightarrow \quad & \displaystyle \max_{t \in \mathbb{R}} \Bigg\{\varepsilon N t - \sum_{i \in [N]} \Bigg(t +  \max\Bigg\{\frac{w_1(\hat{\xi}_{i,1}-x_1)}{\max\{w_1, w_2\}}, \frac{w_2(\hat{\xi}_{i,2}-x_2)}{\max\{w_1, w_2\}}\Bigg\}\Bigg)^+ \Bigg\} \geq \theta N \nonumber \\
\Longleftrightarrow \quad & \displaystyle \max_{t \in \mathbb{R}} \Bigg\{\varepsilon N t - \sum_{i \in [N]} \Bigg(t- \min\Bigg\{\frac{w_1(x_1-\hat{\xi}_{i,1})}{\max\{w_1, w_2\}}, \frac{w_2(x_2-\hat{\xi}_{i,2})}{\max\{w_1, w_2\}}\Bigg\}\Bigg)^+ \Bigg\} \geq \theta N. \label{eq:cvar_ref_3}
\end{align}
The objective function of the embedded maximization problem in~\eqref{eq:cvar_ref_3} is piecewise affine and concave in $t$. Consider the permutation $\bm{\pi}(\bm{x})$ that orders the data points $\bmh{\xi}_i$, $ i \in [N]$, such that
$$
\min \Bigg\{\dfrac{w_1(x_1-\hat{\xi}_{i,1})}{\max\{w_1,w_2\}}, \dfrac{w_2(x_2-\hat{\xi}_{i,2})}{\max\{w_1,w_2\}}\Bigg\} \leq \min \Bigg\{\dfrac{w_1(x_1-\hat{\xi}_{j,1})}{\max\{w_1,w_2\}}, \dfrac{w_2(x_2 - \hat{\xi}_{j,2})}{\max\{w_1,w_2\}}\Bigg\} ~~~1 \leq i \leq j \leq N.
$$
The first-order optimality condition for non-smooth optimization then implies that the maximum on the left-hand side of~\eqref{eq:cvar_ref_3} is attained by 
$$
t^\star = \min \Bigg\{\dfrac{w_1(x_1 - \hat{\xi}_{\pi_{\lfloor \varepsilon N \rfloor + 1}(\bm{x}), 1})}{\max\{w_1, w_2\}}, \dfrac{w_2(x_2 - \hat{\xi}_{\pi_{\lfloor \varepsilon N \rfloor + 1}(\bm{x}), 2})}{\max\{w_1, w_2\}}\Bigg\}.
$$
This implies that the worst-case CVaR constraint~\eqref{eq:cvar_ref_3} holds if and only if
\begin{equation}\label{eq:cvar_almost_done}
\displaystyle \sum_{i = 1}^{\varepsilon N}\min \Bigg\{\dfrac{w_1(x_1 - \hat{\xi}_{\pi_i(\bm{x}),1})}{\max \{w_1, w_2\}}, \dfrac{w_2(x_2 - \hat{\xi}_{\pi_i(\bm{x}),2})}{\max \{w_1, w_2\}}\Bigg\} \geq \theta N.
\end{equation}
Note that $\frac{w_1}{\max\{w_1, w_2\}} \leq 1$ in the first term inside the minimum. Hence, a necessary condition for the inequality~\eqref{eq:cvar_almost_done} to hold for any scaling factors $(w_1, w_2)$ is that $\sum_{i = 1}^{\varepsilon N} (x_1 - \hat{\xi}_{\pi_i (\bm{x}),1}) \geq \theta N$; otherwise, the sum of the first terms inside the minima is smaller than $\theta N$. Note that for any permutation $\bm{\pi} (\bm{x})$, the strong law of large numbers implies that $\frac{1}{N}\sum_{i = 1}^{\varepsilon N} \hat{\xi}_{\pi_i (\bm{x}),1}$ converges to a number smaller than or equal to $\rho$ w.p.~1 as $N$ approaches infinity. Since $\frac{1}{N}\sum_{i = 1}^{\varepsilon N} x_1 = \varepsilon x_1$, we thus conclude that $\frac{1}{N}\sum_{i = 1}^{\varepsilon N} (x_1 - \hat{\xi}_{\pi_i (\bm{x}),1})$ converges to a number not exceeding $\varepsilon x_1 - \rho$ w.p.~1 as $N$ approaches infinity. Since $\overline{x}_1 \varepsilon < \rho$ by assumption, this implies that the inequality~\eqref{eq:cvar_almost_done} is violated for all $x_1 \in [\underline{x}_1, \overline{x}_1]$ w.p.~1 as $N$ approaches infinity.
\hfill \Halmos
\vspace{2mm}

\noindent \emph{Proof of \Cref{prop:CVaR better than bonferroni}.}
We proceed in three steps. We first derive the optimal value of the classical chance constrained program associated with~\eqref{data-driven_2} under the true data-generating distribution $\mathbb{P}_0$ (Step~1). This value serves as a lower bound on the optimal value of problem~\eqref{data-driven_2}. We then show that the worst-case CVaR approximation achieves this bound w.p.~1 as $N \to \infty$~(Step~2), whereas the Bonferroni approximation becomes infeasible~(Step~3).

{\em Step~1.} Since $\rho < \varepsilon$, a similar argument as in the proof of \Cref{prop:bonferroni better than CVaR} allows us to conclude that the optimal value of the classical chance constrained program under the true data-generating distribution $\mathbb{P}_0$ is $\underline{x}$, which is attained by the solution $(x_1, x_2, x_3) = (\underline{x}, \underline{x}, \underline{x})$.

{\em Step~2.} By \Cref{prop:min-signed-distances}, the solution $\bm{x} = (x_1, x_2, x_3) = (\underline{x}, \underline{x}, \underline{x})$ is feasible in the worst-case CVaR approximation with scaling factors $(w_1, w_2) = (\frac{1}{2}, \frac{1}{2})$ if and only if
\begin{equation}\label{eq:cvar_rule_1}
\dfrac{1}{N} \sum_{i = 1}^{\varepsilon N} \mathbf{min\textbf{-}dist} (\hat{\xi}_{\pi_i (\bm{x})}, \mathcal{H}_1(\bm{x}), \mathcal{H}_2(\bm{x})) \geq \theta,
\end{equation}
where $\mathcal{H}_1 (\bm{x}) = \mathcal{H}_2(\bm{x}) = \{\xi \mid \xi \geq \underline{x}\}$, and the permutation $\bm{\pi} (\bm{x})$ orders the data points $\hat{\xi}_i$ such that $\hat{\xi}_1, \ldots, \hat{\xi}_I = 1$, $I \in [N] \cup \{ 0 \}$, and $\hat{\xi}_{I+1}, \ldots, \hat{\xi}_N = 0$. Since $\mathbf{min\textbf{-}dist} (\hat{\xi}_i, \mathcal{H}_1 (\bm{x}), \mathcal{H}_2 (\bm{x})) = \underline{x} - 1$ for $i = 1, \ldots, I$ and $\mathbf{min\textbf{-}dist} (\hat{\xi}_i, \mathcal{H}_1 (\bm{x}), \mathcal{H}_2 (\bm{x})) = \underline{x}$ for $i = I + 1, \ldots, N$,~\eqref{eq:cvar_rule_1} holds if and only if
\begin{equation*}
\dfrac{1}{N}\left(\min\{\varepsilon N, I\} (\underline{x} - 1) + (\varepsilon N - I)^+ \underline{x} \right) \geq \theta.
\end{equation*}
Note that $I / N \to \rho$ w.p.~1 as $N \to \infty$ by the strong law of large numbers. Since $\rho < \varepsilon$ and $\theta \to 0$ as $N \to \infty$, the above inequality is thus satisfied w.p.~1 as $N \to \infty$ as long as $\rho (\underline{x} - 1) + (\varepsilon - \rho) \underline{x} = \varepsilon \underline{x} - \rho$ is strictly positive. This is the case since $\rho < \underline{x} \varepsilon$ by assumption.

{\em Step~3.} Observe that the Bonferroni approximation is infeasible if $\varepsilon_1 \leq I/N$ because the first individual chance constraint $\mathbb{P} [x_1 > \tilde{\xi}] \geq 1-\varepsilon_1 ~\forall \mathbb{P} \in \mathcal{F}(\theta)$ is already violated under the empirical distribution. For the same reason, the Bonferroni approximation is infeasible if $\varepsilon_2 \leq I/N$. We next show that when $N \to \infty$, any pair of Bonferroni weights $(\varepsilon_1, \varepsilon_2)$ satisfying $\varepsilon_1 + \varepsilon_2 = \varepsilon$ also satisfies $\min\{\varepsilon_1, \varepsilon_2\} \leq I/N$ w.p.~1, that is, at least one of the two individual chance constraints is violated. Indeed, we have $\min\{\varepsilon_1, \varepsilon_2\} \leq \varepsilon/2$ and $\rho > \varepsilon/2$ by assumption, and $I/N \to \rho$ w.p.~1 as $N \to \infty$ by the strong law of large numbers.
\hfill \Halmos
\vspace{2mm}

\noindent \emph{Proof of \Cref{prop:ALSO-X better than CVaR}.} $\;$ 
We proceed in three steps. We first derive the optimal value of the classical chance constrained program associated with~\eqref{data-driven_3} under the true data-generating distribution $\mathbb{P}_0$ (Step~1). Since $\mathcal{F} (\theta)$ contains $\mathbb{P}_0$ w.p.~1 as $N \to \infty$, this value bounds the optimal value of problem~\eqref{data-driven_3} from below. We then show that the ALSO-X approximation achieves this bound w.p.~$1$ as $N \to \infty$~(Step~2), whereas the worst-case CVaR approximation becomes infeasible~(Step~3).

{\em Step~1.} Under $\mathbb{P}_0$, we have $0 > \tilde{\xi}_1$ and $0 < \tilde{\xi}_2$ jointly with probability $1 - \rho$, and $x = 0$ therefore satisfies the chance constraint under $\mathbb{P}_0$ since $\rho < \varepsilon$. The optimal value of problem~\eqref{data-driven_3} is thus bounded from below by $0$ w.p.~1 as $N \to \infty$.

{\em Step~2.} We prove the statement by \emph{(i)} showing that for any $\delta \geq 0$, the subproblem of the ALSO-X approximation,
\[
\begin{array}{cll}
\displaystyle \min_{x} & \displaystyle \sup_{\mathbb{P} \in \mathcal{F}(\theta)} \mathbb{E}_{\mathbb{P}}[(\max\{\tilde{\xi}_1 - x, ~x - \tilde{\xi}_2\})^+] \\
{\rm s.t.} & 0 \leq x \leq 1, \; x \leq \delta,
\end{array}
\]
is optimized by $x = 0$ w.p.~1 as $N \to \infty$ and \emph{(ii)} verifying that $x = 0$ also satisfies the distributionally robust chance constraint in the ALSO-X approximation for any value of $\delta$ w.p.~1 as $N \to \infty$.

In view of \emph{(i)}, it follows from Theorem~6.3 of \cite{Esfahani_Kuhn_2017} that the objective function of the subproblem evaluates to
\[
\sup_{\mathbb{P} \in \mathcal{F}(\theta)} \mathbb{E}_{\mathbb{P}}[(\max\{\tilde{\xi}_1 - x, ~x - \tilde{\xi}_2\})^+] 
\;=\; \dfrac{1}{N}\sum_{i \in [N]} (\max\{\hat{\xi}_{i,1} - x, ~x - \hat{\xi}_{i,2}\})^+ + \theta.
\]
Fix any $x \in [0,1]$. For $\bmh{\xi}_i = (1, -1)^\top$, we have that $(\max\{\hat{\xi}_{i,1} - x, ~x - \hat{\xi}_{i,2}\})^+ = \max\{1-x, ~x+1\} = x+1$ strictly increases in $x$, whereas for $\bmh{\xi}_i = (-1,1)^\top$, $(\max\{\hat{\xi}_{i,1} - x, ~x - \hat{\xi}_{i,2}\})^+ = 0$ is non-decreasing in $x$. Thus, as long as there exists at least one data point $\bmh{\xi}_i$ such that $\bmh{\xi}_i = (1,-1)$, the objective function of the subproblem strictly increases in $x$, and the subproblem has the unique optimal solution $x^\star(\delta) = 0$. Since $\mathbb{P}_0 [\tilde{\bm{\xi}} = (1, -1)] = \rho > 0$, this happens w.p.~1 as $N \to \infty$.

As for \emph{(ii)}, \Cref{thm:cc equivalent} implies that
\begin{align}
& \displaystyle \mathbb{P} [0 > \tilde{\xi}_1, ~0 < \tilde{\xi}_2] \geq 1-\varepsilon ~~~\forall \mathbb{P} \in \mathcal{F}(\theta) \nonumber \\
\Longleftrightarrow \quad
& \displaystyle \dfrac{1}{N}\sum_{i = 1}^{\varepsilon N} \mathbf{dist}(\bmh{\xi}_{\pi_i(0)}, \bar{\mathcal{S}}(0)) \ge \theta \nonumber \\
\Longleftrightarrow \quad
& \displaystyle \dfrac{1}{N} \big(\min\{\varepsilon N, I\} \cdot 0 + (\varepsilon N - I)^+ \cdot 1\big) = (\varepsilon - I/N)^+\geq \theta, \label{equ:ALSO-X better CVaR}
\end{align}
where the unsafe set satisfies $\bar{\mathcal{S}}(0) = \{\bm{\xi} \mid \xi_1 \geq 0 ~{\rm or}~ \xi_2 \leq 0\}$, where $\bm{\pi}(0)$ is a permutation that orders $\{\bmh{\xi}_i\}_{i \in [N]}$ by their distances to $\bar{\mathcal{S}}(0)$, and where $I$ denotes the number of samples $\bmh{\xi}_i$, $i \in [N]$, that satisfy $\bmh{\xi}_i = (1, -1)^\top$. Here, the last row follows from the fact that $\mathbf{dist}(\bmh{\xi}_i, \bar{\mathcal{S}}(0)) = 0$ if $\bmh{\xi}_i = (1, -1)^\top$ and $\mathbf{dist}(\bmh{\xi}_i, \bar{\mathcal{S}}(0)) = 1$ if $\bmh{\xi}_i = (-1,1)^\top$. Since $I/N \to \rho$ w.p.~1 as $N \to \infty$ by the strong law of large numbers,~\eqref{equ:ALSO-X better CVaR} is satisfied w.p.~1 as $N \to \infty$ since $\rho < \varepsilon$. In summary, $x = 0$ satisfies the distributionally robust chance constraint in the ALSO-X approximation for any value of $\delta$ w.p.~1 as $N \to \infty$, that is, the ALSO-X approximation is asymptotically exact.

{\em Step~3.} 
Because $x = 1$ violates the chance constraint under $\mathbb{P}_0$, it is infeasible in the worst-case CVaR approximation of~\eqref{data-driven_3} w.p.~1 as $N \to \infty$. We thus restrict our attention to $x \in [0,1)$. Using similar arguments as in the proof of \Cref{prop:bonferroni better than CVaR}, one can show that the worst-case CVaR approximation is satisfied for a fixed decision $x \in [0,1)$ if and only if
\begin{align}
& \sum_{i = 1}^{\varepsilon N}\min \Bigg\{\dfrac{w(x - \hat{\xi}_{\pi_i(x),1})}{\max \{w, 1-w\}}, \dfrac{(1-w)(\hat{\xi}_{\pi_i(x),2} - x)}{\max \{w, 1-w\}}\Bigg\} \geq \theta N \nonumber \\
\Longleftrightarrow \quad
& \dfrac{1}{N} \sum_{i = 1}^{\varepsilon N}\min \big\{w(x - \hat{\xi}_{\pi_i(x),1}), (1-w)(\hat{\xi}_{\pi_i(x),2} - x)\big\} \geq \theta  \cdot \max \{w, 1-w\}. \label{eq:cvar_almost_done_1}
\end{align}
We claim that the permutation $\bm{\pi} (x)$ of the data points implies that there is $I \in \mathbb{N}_0$ such that $\hat{\bm{\xi}}_{\pi_i (x)} = (1, -1)$ for $i \in \{ 1, \ldots, I \}$ and $\hat{\bm{\xi}}_{\pi_j (x)} = (-1, 1)$ for $j \in \{ I + 1, \ldots, N \}$. Indeed, we have
\[
\min \big\{w(x - \hat{\xi}_{i,1}), (1-w)(\hat{\xi}_{i,2} - x)\big\} < 0 < \min \big\{w(x - \hat{\xi}_{j,1}), (1-w)(\hat{\xi}_{j,2} - x)\big\}
\]
whenever $\hat{\bm{\xi}}_i = (1, -1)$ and $\hat{\bm{\xi}}_j = (-1, 1)$ since $x - \hat{\xi}_{i,1} = x - 1 < 0 < x + 1 = x - \hat{\xi}_{j,1}$ and $\hat{\xi}_{i,2} - x = -1 - x < 0 < 1 - x = \hat{\xi}_{j,2} - x$.
Hence, the left-hand side of~\eqref{eq:cvar_almost_done_1} can be re-expressed as
\[
\begin{array}{cl}
& \displaystyle \dfrac{1}{N} \sum_{i = 1}^{\varepsilon N}\min \big\{w(x - \hat{\xi}_{\pi_i(x),1}), (1-w)(\hat{\xi}_{\pi_i(x),2} - x)\big\} \\[3mm]
= & \displaystyle \min\{\varepsilon, I/N\} \cdot \min \big\{w(x - 1), (1-w)(-1 - x)\big\} + (\varepsilon - I/N)^+ \cdot \min \big\{w(x + 1), (1-w)(1 - x)\big\} \\[3mm]
\triangleq & g(w, x).
\end{array}
\]
Thus for any fixed $x \in [0,1)$, the worst-case CVaR approximation is feasible if and only if $g(w, x) \geq \theta \cdot \max \{w, 1-w\}$ for some $w \in (0,1)$. Since
$w(x-1) \leq (1-w)(-1-x)$ if and only if $w \geq \frac{1+x}{2}$ whereas $w(x+1) \leq (1-w)(1-x)$ if and only if $w \leq \frac{1-x}{2}$, for any fixed $x \in [0,1)$ the univariate function $g(\cdot,x)$
has up to two break points $0 < \frac{1-x}{2} \leq \frac{1+x}{2} < 1$.
Recall that as $N \to \infty$, $I/N \to \rho < \varepsilon$ w.p.~1. Consequently, if 
$0 < w \leq \frac{1-x}{2}$, then 
\[
\begin{array}{rl}
g(w,x) = & \min\{\varepsilon, I/N\}(1-w)(-1-x) + (\varepsilon - I/N)^+w(x+1) \\
\xrightarrow{N \to \infty} & \varepsilon(x+1) \cdot w - \rho(x+1) \text{~~~w.p.~1},
\end{array}
\]
which is non-decreasing in $w$ for fixed $x$ because $x+1 \geq 0$; if $\frac{1-x}{2} \leq w \leq \frac{1+x}{2}$, then
\[
\begin{array}{rl}
g(w,x) = & \min\{\varepsilon, I/N\}(1-w)(-1-x) + (\varepsilon - I/N)^+(1-w)(1-x) \\
\xrightarrow{N \to \infty} & (\rho(x+1) - (\varepsilon - \rho)(1-x)) \cdot w + \varepsilon - \varepsilon x - 2\rho \text{~~~w.p.~1},
\end{array}
\]
which is non-decreasing in $w$ for fixed $x$ because $\rho > \varepsilon - \rho$ and $x+1 \geq 1-x$; if $\frac{1+x}{2} \leq w < 1$, then 
\[
\begin{array}{rl}
g(w,x) = & \min\{\varepsilon, I/N\}w(x-1) + (\varepsilon - I/N)^+(1-w)(1-x) \\
\xrightarrow{N \to \infty} & \varepsilon (x-1) \cdot w - (\varepsilon - \rho)(x-1) \text{~~~w.p.~1},
\end{array}
\]
which is non-increasing in $w$ for fixed $x$ as $x-1 \leq 0$. Thus, for any fixed $x \in [0,1)$ we have
\[
\max_{w \in (0,1)} \; g(w,x) = g\left(\frac{1+x}{2},x\right) = (1-x)\left(\frac{(1-x)\varepsilon}{2}-\rho\right) \text{~~~w.p.~1 as $N \to \infty$}, 
\]
which implies that
\[
\max_{x \in [0,1)} \; \max_{w \in (0,1)} \; g(w,x) = \max_{x \in [0,1)} \; (1-x)\left(\frac{(1-x)\varepsilon}{2}-\rho\right) < 0 \text{~~~w.p.~1 as $N \to \infty$},
\]
where the inequality follows from the fact that for any fixed $x \in [0,1)$, $(1-x) > 0$ and $\frac{(1-x)\varepsilon}{2} - \rho \leq \frac{\varepsilon}{2} - \rho < 0$. Thus, for any fixed $x \in [0,1)$ the inequality~\eqref{eq:cvar_almost_done_1}, and therefore the worst-case CVaR approximation to~\eqref{data-driven_3}, is violated for any choice of scaling factors w.p.~1 as $N \to \infty$
\hfill \Halmos
\vspace{2mm}

\noindent \emph{Proof of \Cref{prop:CVaR better than ALSO-X}.} $\;$ 
We proceed in three steps. We first derive the optimal value of the classical chance constrained program associated with~\eqref{data-driven_4} under the true data-generating distribution $\mathbb{P}_0$ (Step~1). This value serves as a lower bound on the optimal value of problem~\eqref{data-driven_4}. We then show that the worst-case CVaR approximation achieves this bound w.p.~1 as $N \to \infty$~(Step~2), whereas the ALSO-X approximation becomes infeasible~(Step~3).

{\em Step~1.} Since $\rho < \varepsilon$, $x \in [0,1)$ is feasible under $\mathbb{P}_0$. Hence, the optimal value of problem~\eqref{data-driven_3}, although not attained, is bounded from below by $-1$ w.p.~1 as $N \to \infty$.

{\em Step~2.} The choice $(w_1, w_2) = (\frac{1}{2}, \frac{1}{2})$ in the worst-case CVaR approximation to~\eqref{data-driven_4} coincides with the choice $\bm{w}^\star$ defined in \Cref{prop:joint cc worst-case CVaR approximation_1}, and \Cref{prop:min-signed-distances} implies that a solution $x\in [0,1)$ is feasible in the worst-case CVaR approximation if and only if
\begin{equation}\label{eq:cvar_rule_4}
\dfrac{1}{N} \sum_{i = 1}^{\varepsilon N} \mathbf{min\textbf{-}dist}(\bmh{\xi}_{\pi_i(x)}, \{ \mathcal{H}_1(x), \mathcal{H}_2(x) \}) \geq \theta,
\end{equation}
where $\mathcal{H}_1(x) = \{\bm{\xi} \mid \xi_1 \geq x\}$, $\mathcal{H}_2(x) = \{\bm{\xi} \mid \xi_2 \leq x\}$, and $\bm{\pi}(x)$ orders the data points $\bmh{\xi}_i$ by their minimum signed distances to $\{ \mathcal{H}_1(x), \mathcal{H}_2(x) \}$. Note that for $\bmh{\xi}_i = (1,1)^\top$, we have $\mathbf{sgn\textbf{-}dist} (\bmh{\xi}_i, \mathcal{H}_1(x)) = x-1$ and $\mathbf{sgn\textbf{-}dist} (\bmh{\xi}_i, \mathcal{H}_2(x)) = 1 - x$; that is, $\mathbf{min\textbf{-}dist} (\bmh{\xi}_i, \mathcal{H}_1(x), \mathcal{H}_2(x)) = x-1 < 0$. Likewise, for $\bmh{\xi}_j = (-1,1)^\top$, we have $\mathbf{sgn\textbf{-}dist} (\bmh{\xi}_j, \mathcal{H}_1(x)) = x + 1$ and $\mathbf{sgn\textbf{-}dist} (\bmh{\xi}_j, \mathcal{H}_2(x)) = 1 - x$; that is, $\mathbf{min\textbf{-}dist} (\bmh{\xi}_j, \mathcal{H}_1(0), \mathcal{H}_2(0)) = 1 - x > 0$. In other words, for any $i \in \{1,\dots, I\}$ and $j \in \{I+1,\dots,N\}$, $x-1 = \mathbf{min\textbf{-}dist} (\bmh{\xi}_i, \mathcal{H}_1(x), \mathcal{H}_2(x)) < 0 < \mathbf{min\textbf{-}dist} (\bmh{\xi}_j, \mathcal{H}_1(x), \mathcal{H}_2(x)) = 1 - x$, and therefore $\bm{\pi}(x)$ orders the data points $\bmh{\xi}_i$ such that $\bmh{\xi}_1, \ldots, \bmh{\xi}_I = (1,1)^\top$, $I \in [N] \cup \{0\}$, and $\bmh{\xi}_{I+1}, \ldots, \hat{\xi}_N = (-1,1)^\top$. Hence, the inequality~\eqref{eq:cvar_rule_4} holds if and only if $\min\{\varepsilon, I/N\} \cdot (x-1) + (\varepsilon - I/N)^+ \cdot (1-x) \geq \theta$. Note that by the strong law of large numbers, $I/N \to \rho$ w.p.~1 as $N \to \infty$. Since $\rho < \varepsilon/2$ and $\theta \to 0$ as $N \to \infty$, the above inequality becomes $\rho(x-1)+(\varepsilon-\rho)(1-x) = (\varepsilon-2\rho)(1-x) \geq \theta$
and is thus satisfied by any $x \in [0,1)$ w.p.~1 as $N \to \infty$.

{\em Step~3.} For any $\delta \geq -1$, the ALSO-X approximation solves the subproblem
\[
\begin{array}{cll}
\displaystyle \min_{x} & \displaystyle \sup_{\mathbb{P} \in \mathcal{F}(\theta)} \mathbb{E}_{\mathbb{P}}[(\max\{\tilde{\xi}_1 - x, x - \tilde{\xi}_2\})^+] \\
{\rm s.t.} & 0 \leq x \leq 1, \; -x \leq \delta.
\end{array}
\]
As in the proof of \Cref{prop:ALSO-X better than CVaR}, the objective function of the subproblem satisfies
\[
\sup_{\mathbb{P} \in \mathcal{F}(\theta)} \mathbb{E}_{\mathbb{P}}[(\max\{\tilde{\xi}_1 - x, x - \tilde{\xi}_2\})^+] = \dfrac{1}{N}\sum_{i \in [N]} (\max\{\hat{\xi}_{i,1} - x, x - \hat{\xi}_{i,2}\})^+ + \theta.
\]
For any $x \in [0,1]$, 
if $\bmh{\xi}_i = (1,1)^\top$, then $(\max\{\hat{\xi}_1 - x, x - \hat{\xi}_2\})^+ = 1-x$ strictly decreases in $x$; if $\bmh{\xi}_i = (-1,1)^\top$, then $(\max\{\hat{\xi}_1 - x, x - \hat{\xi}_2\})^+ = 0$ is non-increasing in $x$. Thus, as long as there exists at least one data point $\bmh{\xi}_i$ such that $\bmh{\xi}_i = (1,1)^\top$, then the objective function of the subproblem strictly decreases in $x$, and the subproblem always returns the unique solution $x^\star(\delta) = 1$. Since $\mathbb{P}_0 [\tilde{\bm{\xi}} = (1, 1)] = \rho > 0$, this happens w.p.~1 as $N \to \infty$. Because $x = 1$ is already infeasible to the classical chance constrained program under $\mathbb{P}_0$, with probability going to $1$ as $N \to \infty$, it cannot be feasible in the ALSO-X approximation of the distributionally robust chance constrained program~\eqref{data-driven_4}. Hence, we have to increase $\delta$ to infinity, which implies that the ALSO-X approximation fails to find a feasible solution.
\hfill \Halmos
\vspace{2mm}

\noindent \emph{Proof of \Cref{prop:ALSO-X better than Bonferroni}.} $\;$ Let $I$ denote the number of samples $\bmh{\xi}_i$, $i \in [N]$, that satisfy $\bmh{\xi}_i = (1,-1)^\top$. The Bonferroni approximation is infeasible if $\varepsilon_1 \leq I/N$ (resp., $\varepsilon_2 \leq I/N$) because the first (resp., second) individual chance constraint $\mathbb{P} [x > \tilde{\xi}_1] \geq 1-\varepsilon_1 ~\forall \mathbb{P} \in \mathcal{F}(\theta)$ (resp., $\mathbb{P} [x < \tilde{\xi}_2] \geq 1-\varepsilon_2 ~\forall \mathbb{P} \in \mathcal{F}(\theta)$) is already violated under the empirical distribution. A similar argument as in the proof of \Cref{prop:CVaR better than bonferroni} thus allows us to conclude that when $N \to \infty$, at least one of the two individual chance constraints is violated w.p.~1 for any pair of $(\varepsilon_1, \varepsilon_2)$ satisfying $\varepsilon_1 + \varepsilon_2 = \varepsilon$.
\hfill \Halmos
\vspace{2mm}

\noindent \emph{Proof of \Cref{prop:bonferroni better than ALSO-X}.} $\;$ Fix any $x \in [0,1)$, and denote by $\mathcal{S}_1(x) = \{\bm{\xi} \mid x > \xi_1\}$ and $\mathcal{S}_2(x) = \{\bm{\xi} \mid x < \xi_2\}$ the two safety sets of the Bonferroni approximation. If $\bmh{\xi}_i = (1,1)^\top$, then $\bmh{\xi}_i \in \bar{\mathcal{S}}_1(x) \cap \mathcal{S}_2(x)$ with $\mathbf{dist}(\bmh{\xi}_i, \bar{\mathcal{S}}_1(x)) = 0$ and $\mathbf{dist}(\bmh{\xi}_i, \bar{\mathcal{S}}_2(x)) = 1-x$. Likewise, if $\bmh{\xi}_i = (-1,1)^\top$, then $\bmh{\xi}_i \in \mathcal{S}_1(x) \cap \mathcal{S}_2(x)$ with $\mathbf{dist}(\bmh{\xi}_i, \bar{\mathcal{S}}_1(x)) = 1+x$ and $\mathbf{dist}(\bmh{\xi}_i, \bar{\mathcal{S}}_2(x)) = 1-x$. Under the appropriate permutations $\bm{\pi}^1(x)$ and $\bm{\pi}^2(x)$, \Cref{thm:cc equivalent} then implies that $x$ satisfies both chance constraints of the Bonferroni approximation if and only if
$\frac{1}{N} \sum_{i = 1}^{\varepsilon_1 N}\mathbf{dist}(\bmh{\xi}_{\pi^1_i(x)}, \bar{\mathcal{S}}_1(x)) = (\varepsilon_1 - I/N)^+ (1+x) \geq \theta$ and $\frac{1}{N} \sum_{i = 1}^{\varepsilon_2 N}\mathbf{dist}(\bmh{\xi}_{\pi^2_i(x)}, \bar{\mathcal{S}}_2(x)) = \varepsilon_2 (1-x) \geq \theta$, where $I$ denotes the number of samples $\bmh{\xi}_i$, $i \in [N]$, that satisfy $\bmh{\xi}_i = (1,1)^\top$. Since $x \in [0,1)$, both constraints are satisfied under any choice $\varepsilon_1 \in (\rho, \varepsilon)$ and $\varepsilon_2 = \varepsilon - \varepsilon_1$ w.p.~1 as $N \to \infty$.
\hfill \Halmos
\vspace{2mm}

\end{appendices}

%%%%%%%%%%%%%%%%%
%%%%%%%%%%%%%%%%%
\end{document}